\documentclass[preprint,12pt,authoryear]{elsarticle}

\usepackage{microtype}

\usepackage{tikz}
\usetikzlibrary{arrows.meta,positioning,shadows,fit,shapes.symbols,shapes.geometric}

\usepackage{multirow}
\usepackage{graphicx}
\usepackage{subfigure}
\usepackage{makecell}
\usepackage{amsmath}
\usepackage{url}
\usepackage{hyperref}
\usepackage{enumitem} 
\usepackage{float} 

\usepackage{amssymb}
\usepackage{amsmath}

\usepackage[ruled,vlined]{algorithm2e}

\usepackage{xcolor}

\journal{Computational and Applied Mathematics}

\begin{document}

\begin{frontmatter}



\title{A Systematic Framework for Evaluating Topological Representations in Single-Cell Classification}



\author[1,2]{Rocío Picón-González \corref{cor1}}

\author[1,2]{Salvador Chulián}

\author[1,3]{Ana Niño-López}

\author[4]{Álvaro Martínez-Rubio}

\author[1,2]{María Rosa Durán} 

\address[1]{Biomedical Research and Innovation Institute of Cádiz (INiBICA), Puerta del Mar University Hospital, Cádiz 11009, Spain}

\address[2]{Department of Mathematics, Universidad de Cádiz, Puerto Real 11510, Spain}

\address[3]{Pediatric Health Research Institute Niño Jesús University Children's Hospital (IPIS-NJ), Madrid, Spain}

\address[4]{Institut Curie, PSL University, Paris, France}

\cortext[cor1]{Corresponding author: \texttt{rocio.picon@uca.es}
}

\begin{abstract}
Recent advances in biomedicine generate high-dimensional single-cell data that describe cellular heterogeneity with unprecedented detail, but their geometric complexity and non-linear structure often limit the effectiveness of conventional statistical tools. Topological Data Analysis (TDA) provides a mathematical framework for characterizing the shape of data through persistent homology, which extracts structural features such as connected components and cycles across multiple scales. In this work, we propose a systematic two-level framework for evaluating topological representations in high-dimensional single-cell classification. The first level (\(R_1\)) performs statistical screening of topological descriptors based on separability between clinical groups, whereas the second level (\(R_2\)) evaluates their predictive utility in supervised classification models. This design makes it possible to compare representations not only in terms of discriminative performance, but also in terms of robustness to analytical choices. We illustrate the framework using bone marrow flow cytometry data from pediatric acute lymphoblastic leukemia, with a particular focus on relapse stratification. The results show that different topological representations vary substantially in both statistical separability and predictive stability, with Betti Curves and Persistence Silhouettes showing more robust behavior than Persistence Images in this cohort. Overall, the study provides a reproducible methodological framework for the systematic comparison of topological descriptors in complex biomedical point clouds.
\end{abstract}

\begin{graphicalabstract}
\includegraphics[width=1\textwidth]{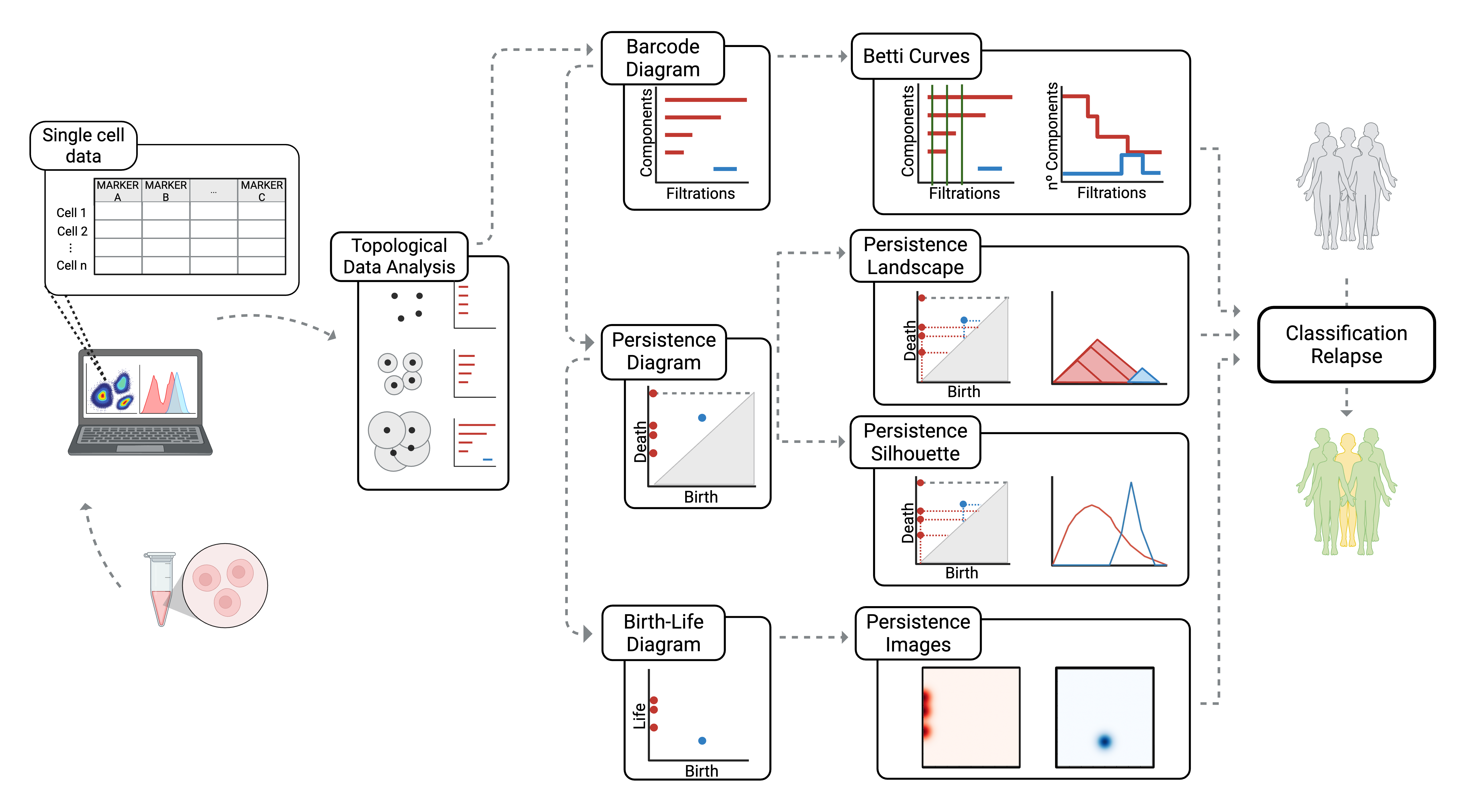}
\end{graphicalabstract}

\begin{highlights}
\item A two-level framework is proposed to evaluate topological representations in single-cell classification.
\item Persistent homology descriptors are compared across representations and parameters.
\item Betti Curves and Persistence Silhouettes show more stable behavior than Persistence Images..
\item Combining \(F_2\), AUC and confusion matrices improves model assessment.
\item Topological descriptors show heterogeneous behavior across analytical settings.
\end{highlights}


\end{frontmatter}



In recent years, biomedicine has undergone a transformation driven by technologies capable of generating high-dimensional, high-resolution data at the single-cell level. This data allows for an exhaustive description of cellular heterogeneity, opening new avenues for precision medicine \citep{wolde2025current}. At the same time, it poses a significant analytical challenge: geometric complexity, the non-linear nature of the data, and the presence of multiscale structures often hinder conventional statistical tools from fully exploiting the available biological information \citep{HernandezLemus2025,applications}.

In this context, Topological Data Analysis (TDA) has emerged as a mathematical framework specifically designed to study the ``shape'' of data. Through tools such as persistent homology, TDA identifies structural properties across scales, such as connected components or cycles. This translates into descriptors capable of capturing global organization and complex relationships within cell populations \citep{EdelsbrunnerHarer2010,ZomorodianCarlsson2005,Chazal2016}. Recent reviews emphasize its growing relevance in single-cell biology, where it has been proposed as a promising pathway to reveal cellular trajectories, transitional states and collective patterns that are not always detectable through linear or purely local approximations \citep{Levenson2024, HernandezLemus2025}.

Several studies have already shown the potential of TDA to extract relevant information from complex biological systems. In oncology, recent reviews have highlighted its growing range of applications, from tissue architecture and tumor organization to clinically relevant pattern discovery \citep{applications}. Specifically, \cite{Stolz2024} proposed a relational persistent homology framework for multispecies data; \cite{Yang2025} employed topological tools to classify tumor-immune interactions; and \cite{BonillaCarpioTrenado2020}  showed that topology can capture complex dynamic patterns  in collective cell motion. In cytometry and related single-cell settings, persistent homology has also been used to identify clinically relevant structures in high-dimensional cellular data \citep{determining}.  In the specific context of acute lymphoblastic leukemia, both \cite{acute} and \cite{bib1} suggest that topological descriptors can contribute to relapse-related classification. Collectively, these works suggest that the latent geometry and topology of biological data contain potentially useful information for characterization and classification tasks.

This scenario also reveals an important methodological challenge. Although TDA has already produced promising results in biomedical research, analytical decisions are not always evaluated in a systematic way. As a result, it is often difficult to determine how much of the observed performance depends on the representation itself and how much depends on specific choices such as homological dimension, parameterization, or classification model. This makes it more difficult to compare studies and to identify which topological descriptors remain stable across different analytical settings.

In response, we propose a general methodological framework aimed at studying the impact of these decisions in a structured way. Rather than focusing only on final predictive performance, the framework is designed to evaluate the sensitivity of different topological representations to analytical choices and to compare their robustness under a common workflow. \mbox{Specifically}, we introduce a two-level analytical scheme (\(R_1\) and \(R_2\)) that allows us to examine how parameter choices and variable combinations affect the extracted topological signal and its utility for classification. In this sense, the study seeks not only to compare representations, but also to contribute to a more systematic strategy for applying TDA in complex biomedical data.

As a case study, we apply this framework to bone marrow flow cytometry data from pediatric leukemia. Cytometry represents a paradigmatic example of single-cell data, as each sample consists of large cell populations described by multiple markers, forming high-dimensional point clouds whose internal organization may reflect subtle biological and clinical states. Our work aligns with this research line, but with a particular emphasis on designing a \mbox{generalizable} methodological strategy to compare different data representations and identify the most suitable configurations for topological analysis and supervised classification.
Beyond comparing topological representations, the main methodological contribution of this work lies in separating statistical topological screening from supervised predictive evaluation within a single framework. This makes it possible to assess not only discriminative performance, but also the robustness of each representation to analytical choices.

\section*{Methods}

\subsection*{\textbf{Persistent Homology}}

Let $X \in \mathbb{R}^{N \times M}$ be a preprocessed sample and 
\[
\mathcal{P} = \{x_1,\dots,x_N\} \subset \mathbb{R}^M
\]
be the point cloud associated, where N represents the number of observations and M the number of characteristics observed. We equip $\mathcal{P}$ with the Euclidean metric
$d(x_i,x_j)=\lVert x_i-x_j\rVert_2$.

To characterize the global structure of the point cloud, we represent it as a simplicial complex, a combinatorial structure composed of points (0-simplices), edges (1-simplices) and higher-dimensional simplices that connect observations according to a distance criterion. Since the choice of a specific distance threshold is often arbitrary, we construct a simplicial filtration, which is defined as an ordered sequence of nested complexes \citep{EdelsbrunnerHarer2010, ZomorodianCarlsson2005, Choi2025}. This multi-scale approach allows for the systematic tracking of topological features across the filtration.

Specifically, for a scale parameter $\varepsilon \ge 0$, the Vietoris-Rips complex $\mathrm{VR}(\mathcal{P},\varepsilon)$ is defined as:

\[
\mathrm{VR}(\mathcal{P},\varepsilon)
=
\bigl\{
\sigma \subseteq \mathcal{P} : d(x_i,x_j) \le \varepsilon 
\ \forall x_i,x_j \in \sigma
\bigr\}.
\]

The condition $\varepsilon_1 \le \varepsilon_2$ implies 
$\mathrm{VR}(\mathcal{P},\varepsilon_1) \subseteq \mathrm{VR}(\mathcal{P},\varepsilon_2)$, 
yielding an increasing sequence of complexes that allows for the multiscale tracking of topological features 
\citep{EdelsbrunnerHarer2010, Chazal2016}.

For each dimension $k$, we compute the birth times $b_i$, death times $d_i$ and persistence $p_i=d_i-b_i$ of the homological classes, summarized in a persistence diagram

\[
D_k = \{(b_i,d_i)\}_{i=1}^{n_k},
\]

where \(n_k\) denotes the number of persistence pairs in dimension \(k\).

This study examines dimensions $k=0$ and $k=1$. $H_0$ persistence tracks the connectivity and clustering of the point cloud, while $H_1$ persistence identifies the formation and dissolution of cycles or loops \citep{SkrabaTurner2020}. We also evaluate a combined $H_0{+}H_1$ representation by concatenating the descriptors from both dimensions.

\begin{figure}[H]
    \centering
    \includegraphics[width=0.99\textwidth]{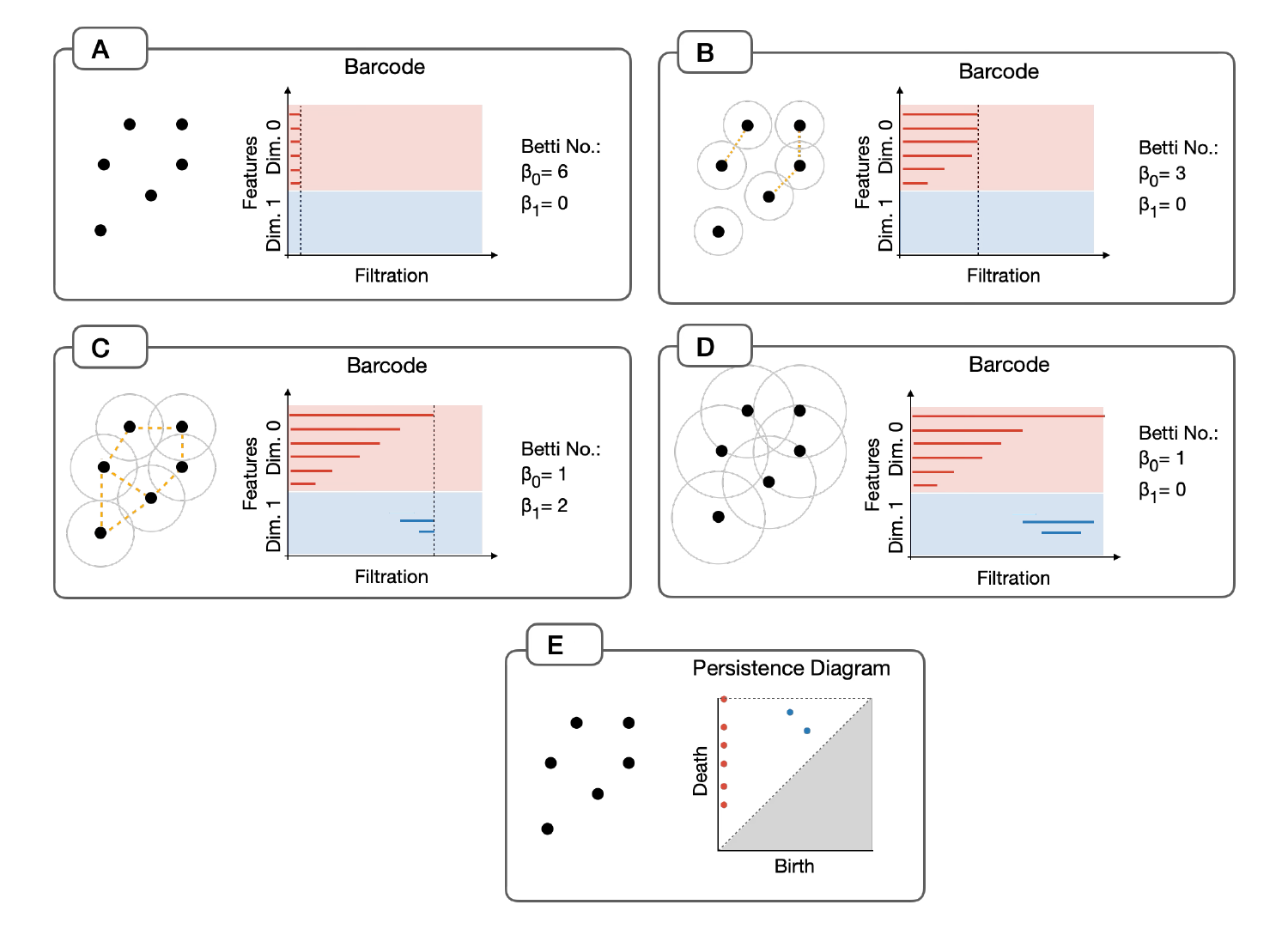}
    \caption{\textbf{Persistent Homology construction via Vietoris-Rips filtration.} \mbox{(\textbf{A-D})} Evolution of the simplicial complex $VR(P,\epsilon)$ as a function of the scale parameter $\epsilon\geq0$. As the neighborhood radius increases, new edges and faces are formed, modifying the topological structure of the system. Structural properties are quantified through the Betti number ($\beta_k$), which represents the rank of the $k$-th homology group $H_k$: $\beta_0$ denotes the number of connected components and $\beta_1$ corresponds to the number of cycles present in the complex for a given value of $\epsilon$. (\textbf{E}) This results in the persistence diagram $D_k$, where each point ($b_i, d_i$) represents the birth and death of a topological feature across the filtration scale.}
    \label{fig:fig_vr}
\end{figure}

\subsection*{\textbf{Topological Vectorizations}}
To integrate multiscale information from persistence diagrams $D_0$ and $D_1$ into statistical and machine learning frameworks, we transform these summaries into four distinct vector representations: Betti Curves, Persistence Landscapes, Persistence Silhouettes and Persistence Images. These methods map topological features into numerical descriptors while preserving the structural properties of the point cloud.

\subsubsection*{Betti Curves}

Given a persistence diagram $D_k = \{(b_i,d_i)\}$, the Betti curve in dimension $k$ is defined as:
\[
\beta_k(t) = \#\{\, (b_i,d_i)\in D_k : b_i \le t < d_i \,\}.
\]
This function counts the number of persistence intervals active in a given filtration step $t$. By discretizing the domain into a sequence $t_1 < \cdots < t_L$, the representation is obtained as a vector.
\[
(\beta_k(t_1),\dots,\beta_k(t_L)) \in \mathbb{R}^L.
\]
Betti curves provide a straightforward global description of how topological features evolve across the filtration \citep{EdelsbrunnerHarer2010}.

\subsubsection*{Persistence Landscapes}

Persistence Landscapes offer a functional representation of persistence diagrams \citep{Bubenik2015}. For each bar $(b_i,d_i)\in D_k$, we define the triangular function:
\[
\lambda_i(t) = \max\bigl(0, \min(t-b_i,\, d_i-t)\bigr).
\]
The $j$-th landscape corresponds to the $j$-th maximum of the set of functions
$\{\lambda_i(t)\}_i$. This approach ensures stability and allows for the use of standard functional data analysis tools.

\subsubsection*{Persistence Silhouettes}

Silhouettes provide an analytical summary based on a weighted average of triangular functions associated with the bars \citep{Chazal2016}. For weights $\omega_i \ge 0$, the silhouettes function is defined as:
\[
\mathrm{Sil}(t)
= \frac{\sum_i \omega_i \max(0,\min(t-b_i,\, d_i-t))}{\sum_i \omega_i}.
\]
This representation summarizes the collective contribution of homological features, assigning higher importance to features with greater persistence.

\subsubsection*{Persistence Images}

Persistence images transform the diagrams into a stable grid-based vector representation \citep{Adams2017}. First, each point $(b_i,d_i)$ is converted into birth-persistence coordinates $(b_i,p_i)$, where $p_i = d_i - b_i$. On this plane, each point is replaced by a weighted Gaussian kernel:

\[
\phi_i(z) = w_i \exp\!\left(-\frac{\lVert z-(b_i,p_i)\rVert^2}{2\sigma^2}\right),
\]
where $w_i$ is a weight function dependent on persistence. The total density,  $\rho(z) = \sum_i \phi_i(z)$, is evaluated over a regular grid of resolution $(n_x,n_y)$, where \(n_x\) and \(n_y\) correspond to the horizontal and vertical resolutions of the persistence image. The resulting image is then vectorized into $\mathbb{R}^{n_x \cdot n_y}$.

 \begin{figure}[H]
            \centering
            \includegraphics[width=0.74\textwidth]{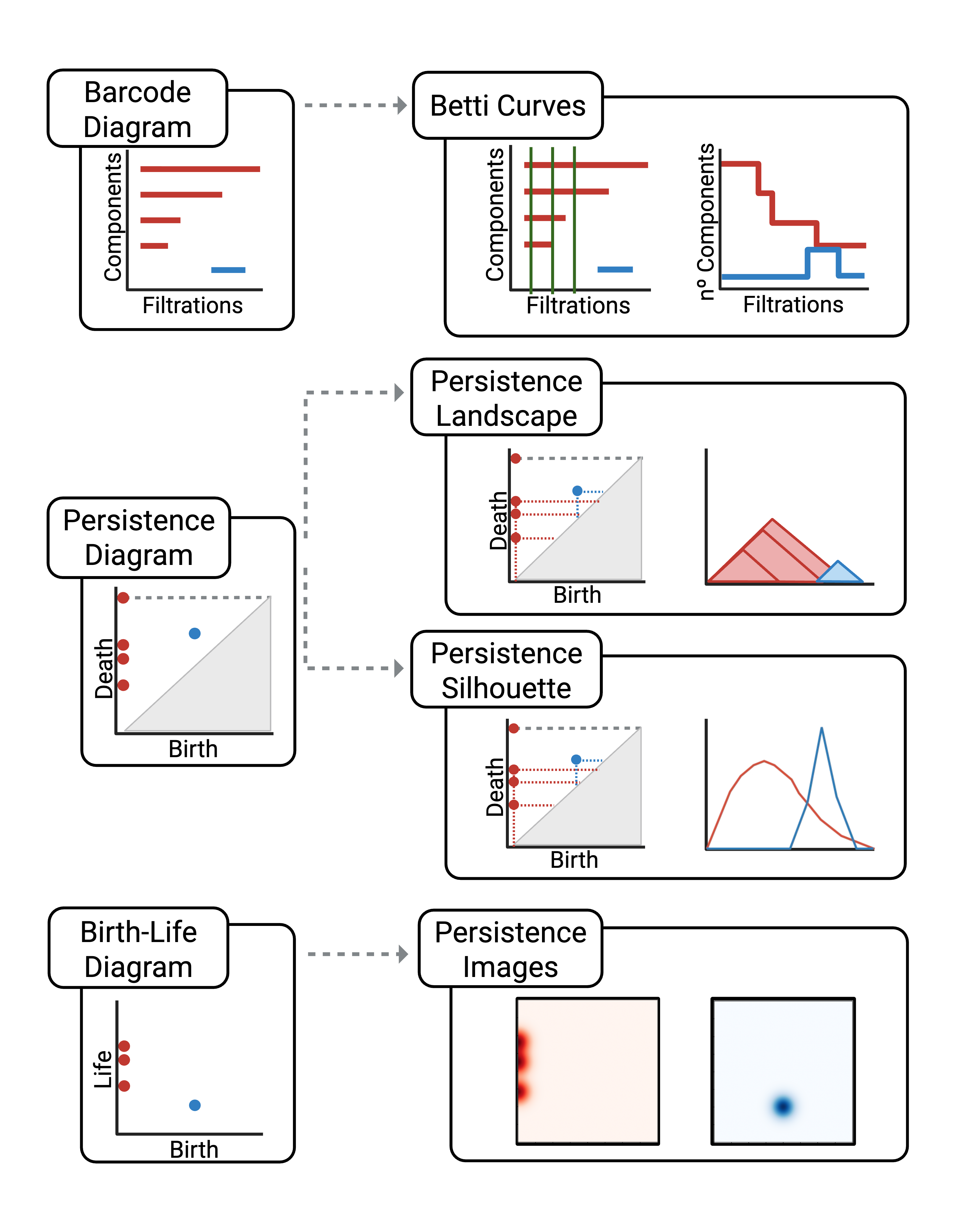}
         \caption{\textbf{Vectorization methods for topological descriptors.} Mapping multiscale information from barcodes and persistence diagrams into functional and vector spaces in $\mathbb{R}^m$, where \(m\) denotes the dimensionality of the resulting representation. The transformations for Betti Curves, Persistence Landscapes, Persistence Silhouettes and Persistence Images are illustrated.}
            \label{fig:repre}
        \end{figure}

To compare the topological representations, we explored both common and representation-specific hyperparameters. Table~\ref{tab:topo_hyperparameters} summarizes the parameters considered in each representation and their corresponding tested values.

\begin{table}[H]
\centering
\small
\newcommand{\hparam}[2]{%
\textbf{#1}\\[-0.2em]
\hspace*{0.8em}%
\begin{minipage}[t]{\dimexpr\linewidth-0.8em\relax}
{\footnotesize\itshape #2}
\end{minipage}%
}
\resizebox{\columnwidth}{!}{%
\begin{tabular}{p{4.6cm} c c c c}
\hline
\textbf{Hyperparameter} 
& \textbf{Betti Curves} 
& \textbf{Pers. Silhouettes} 
& \textbf{Pers. Landscapes} 
& \textbf{Pers. Images} \\
\hline

\hparam{Resolution}{Number of discretization bins or grid cells used to represent the topological summary}
& 100, 125, 150
& 100, 125, 150
& 100, 125, 150
& 10, 25, 50 \\

\hline

\hparam{Smoothing}{Degree of smoothing applied to the functional or vector representation}
& 0, 1, 2
& 0, 1, 2
& 0, 1, 2
& 0, 2 \\

\hline

\hparam{Threshold}{Lower cutoff used to remove low-persistence features}
& 0, p10
& 0, p10
& 0, p10
& 0, p10 \\

\hline

\hparam{Normalization}{Scaling applied to the representation before statistical or supervised analysis}
& none, $L^1$
& none, $L^1$
& none, $L^1$
& none, $L^1$ \\

\hline

\hparam{Landscapes}{Number of landscape functions retained in the descriptor}
& --
& --
& 1, 3, 5
& -- \\

\hline

\hparam{Bandwidth}{Width of the Gaussian kernel used to spread persistence points in the image}
& --
& --
& --
& 0.5, 2 \\

\hline

\hparam{Weights}{Weighting scheme used to modulate the contribution of persistence points}
& --
& --
& --
& const, persist \\

\hline
\end{tabular}
}
\caption{\textbf{Hyperparameter grids explored for each topological representation.} For each hyperparameter, the table reports the values explored for Betti curves, Persistence Silhouettes (Pers. Silhouettes), Persistence Landscapes (Pers. Landscape), and Persistence Images (Pers. Images). In this context, \textit{p10} denotes the 10th percentile of persistence values, $L^1$ indicates normalization by the sum of absolute values, and \textit{const} and \textit{persist} refer to constant and persistence-based weighting schemes, respectively.}
\label{tab:topo_hyperparameters}

\end{table}

\subsection*{\textbf{Mann-Whitney--Fisher Statistical Framework}}

To evaluate the discriminative power of each topological representation without relying on supervised learning models, we implement a bin-by-bin statistical analysis based on the non-parametric Mann-Whitney U test \citep{mann1947test, stat}. 
Let $Z$ be the set of patients, partitioned into two cohorts \(Z_A\) and \(Z_B\).

Given a vector representation $f(X)$ of a sample $X$, each coordinate $j$ (representing a bin or pixel) is analyzed independently. For each coordinate, we compare the distributions:
\[
\{ f_j(X) : X \in Z_{A} \}
\quad\text{and}\quad
\{ f_j(X) : X \in Z_{B} \}
\]
using the Mann-Whitney $U$ statistic. This yields a $p$-value $p_j \in (0,1)$ that quantifies the evidence of separation between both groups at that specific bin.

To obtain a global measure of separability for the entire representation, the individual $p$-values are combined using Fisher's method \citep{fisher}:
\[
T = -2 \sum_{j=1}^{d} \log(p_j),
\]
where $d$ denotes the dimension of the representation. Under the null hypothesis, the statistic $T$ follows a $\chi^2$ distribution with $2d$ degrees of freedom, from which the combined value $p_{\mathrm{Fisher}}$ is derived. 

Finally, we define the separability score $S$ as:
\[
S = -\log_{10}(p_{\mathrm{Fisher}}.
\]
Higher values of $S$ indicate a more consistent statistical separation between the cohorts under comparison (\(Z_A\) and \(Z_B\)). This metric is used within each hyperparameter grid to identify and select the most informative topological configuration.

In this setting, Fisher's method is used as a practical way to aggregate coordinate-wise evidence of separation across the representation. Since the coordinates of a vectorized topological descriptor may be correlated, the resulting score should be interpreted primarily as an aggregate separability measure rather than as a formal inferential test under strict independence assumptions.

\subsection*{\textbf{Classification Models}}
To assess the predictive utility of the selected topological descriptors, we implement a supervised learning stage using two families of classifiers with distinct mathematical foundations: Random Forest and Support Vector Machines.

\subsubsection*{Random Forest}

Random Forest (RF) is an ensemble learning method that constructs a multitude of decision trees during training to mitigate the variance and overfitting characteristic of individual trees \citep{Breiman2001,sas,stat}. The algorithm employs bagging (bootstrap aggregating), where each tree is trained on a random sample of the data drawn with replacement. Feature randomness is introduced by considering only a subset of markers at each node split to ensure diversity among the estimators. As illustrated in the conceptual diagram in Figure~\ref{fig:class}, the final output is obtained through a majority voting scheme that enhances model stability. In our study, we optimized the number of trees and the maximum tree depth as the primary hyperparameters.

\subsubsection*{Support Vector Machine}

\begin{figure}[H]
    \centering
    \includegraphics[width=0.9\textwidth]{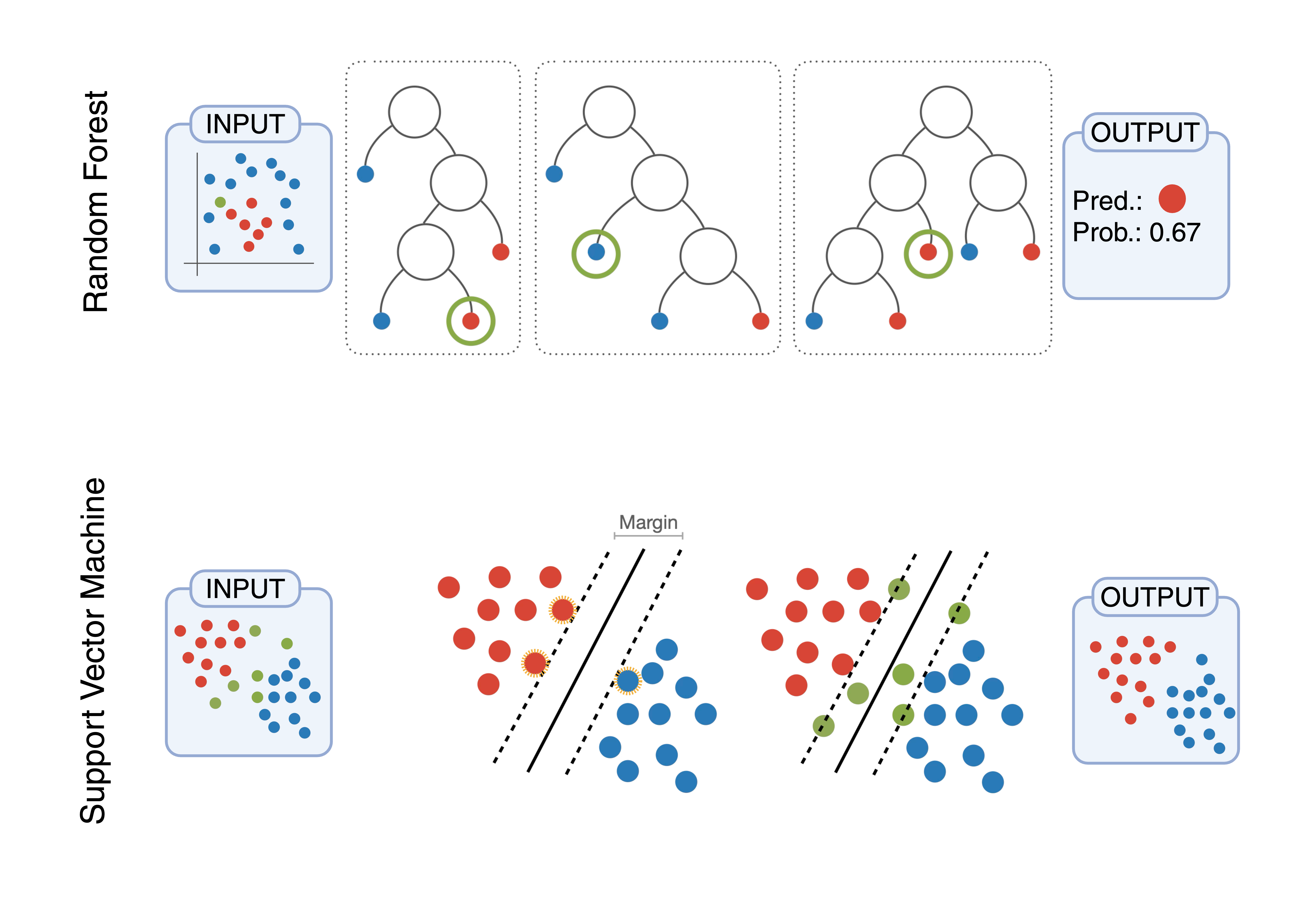}
    \caption{\textbf{Conceptual illustration of Random Forest and Support Vector Machines.} In the Random Forest framework (top), multiple decision trees are trained on random subsets of data and features to generate independent predictions. The final classification is determined by a majority voting scheme. In the Support Vector Machines model (bottom), the solution is decided by the support vectors, which identify the optimal hyperplane that maximizes the margin between classes. New observations, represented in green, are classified according to their relative distance and position with respect to this decision boundary.}
    
    \label{fig:class}
\end{figure}

The Support Vector Machine (SVM) framework identifies an optimal hyperplane that separates the cohorts by maximizing the margin, defined as the distance  between the decision boundary and the nearest observations or support vectors \citep{Cortes1995, prml, stat}. This geometric principle is depicted in Figure \ref{fig:class}, where the support vectors effectively define the solution. Since biological point clouds are rarely linearly separable, we use the kernel trick to map descriptors into higher-dimensional feature spaces \citep{prml, stat}. We evaluated four kernel functions: linear, Radial Basis Function  (RBF), polynomial and sigmoidal. The regularization parameter $C$ was tuned to balance the trade-off between margin maximization and training error minimization.

\subsection*{\textbf{Model Evaluation}}

Predictive performance is assessed through stratified cross-validation (CV), an approach that maintains consistent class distribution across all training and evaluation partitions. This procedure is essential when handling imbalanced datasets, ensuring that each fold accurately reflects the original proportions of the cohorts. By preserving these ratios, the cross-validation process generates a robust estimation of the classifier's capability to generalize across unseen data \citep{Lopez2025Validation, Brishti2025Imbalanced}.

Performance is characterized using standard metrics for binary classification, which are calculated from the predictions obtained in the validation sets \citep{stat,HanleyMcNeil1982,VanRijsbergen1979}. 
We adopt the standard confusion matrix notation: True Positives (TP), True Negatives (TN), False Positives (FP) and False Negatives (FN) \citep{stat}.

\paragraph{Accuracy}
Accuracy (ACC) measures the total proportion of correct predictions and is defined as \cite{stat}:
\[
\mathrm{ACC}=\frac{\mathrm{TP}+\mathrm{TN}}{\mathrm{TP}+\mathrm{TN}+\mathrm{FP}+\mathrm{FN}}.
\]
While intuitive, this metric can be misleading in imbalanced datasets, as it does not distinguish between performance in the majority and minority classes.

\paragraph{Area Under the ROC Curve}
The Area Under the ROC Curve (AUC) quantifies the model's ability to separate both cohorts across all possible decision thresholds\citep{HanleyMcNeil1982}. The ROC curve is defined by the set of pairs $(\mathrm{FPR}(t),\mathrm{TPR}(t))$ as the threshold $t$ varies, with:
\[
\mathrm{TPR}(t)=\frac{\mathrm{TP}(t)}{\mathrm{TP}(t)+\mathrm{FN}(t)}, \qquad
\mathrm{FPR}(t)=\frac{\mathrm{FP}(t)}{\mathrm{FP}(t)+\mathrm{TN}(t)}.
\]
The area under this curve is expressed as:
\[
\mathrm{AUC}=\int_0^1 \mathrm{TPR}(u)\,du,
\]
where $u=\mathrm{FPR}$. This metric provides a global measure of the classifier's discriminative power and is less sensitive to class imbalance than accuracy \citep{HanleyMcNeil1982}. 

\paragraph{$F_2$-measure}
The F-score provides a balance between precision (P) and recall (R), which are defined as \citep{VanRijsbergen1979, Christen2023FMeasure}:
\[
P=\frac{\mathrm{TP}}{\mathrm{TP}+\mathrm{FP}}, \qquad
R=\frac{\mathrm{TP}}{\mathrm{TP}+\mathrm{FN}},
\]
The general family of these scores is defined as:
\[
F_{\beta}=\frac{(1+\beta^2)\,P\,R}{\beta^2 P + R}.
\]
In this study, we use the \(F_2\) variant,
\[
F_2=\frac{5PR}{4P+R}.
\]

This metric assigns four times more weight to recall than to precision, prioritizing the correct identification of relapse cases due to their clinical significance.

\paragraph{Confusion Matrix}
The Confusion Matrix (CM) provides a detailed summary of the classifier's performance \citep{stat}
\[
\begin{pmatrix}
\mathrm{TN} & \mathrm{FP} \\
\mathrm{FN} & \mathrm{TP}
\end{pmatrix},
\]
where the entries correspond to the notation established previously. This representation is fundamental for a direct interpretation of the clinical implications of the models, as it allows for a clear distinction between the types of misclassification errors. 

\subsubsection*{Data Imbalance}

To mitigate the effects of class imbalance, we evaluate class balancing through an oversampling strategy. This procedure involves the replication of existing minority class observations within the training set, a method documented to improve model performance in imbalanced biological datasets \citep{he2009learning}.

We prioritized this replication strategy to preserve the fidelity of the original topological structures within the feature space, avoiding the introduction of synthetic noise that could distort the multiscale relationships captured by the TDA. To maintain the integrity of the evaluation, oversampling is applied strictly to the training partitions of the stratified cross-validation, leaving the validation sets unmodified.

\subsubsection*{Specification Curves}

Specification curves allow the systematic evaluation of how a score varies across the full set of analytical configurations considered in a study \citep{simonsohn2020specification}. Rather than focusing only on the best-performing result, this approach represents the distribution of scores across all valid combinations of methodological choices.

Each configuration is associated with a single summary score, and the full set of configurations is ranked according to that value. The resulting curve allows one to assess whether the observed behavior is concentrated in a narrow subset of highly specific settings or remains relatively stable across a broader region of the analytical space. In this sense, specification curves provide a useful tool for evaluating the sensitivity and robustness of topological descriptors under multiple analytical decisions.

In addition to the ranked performance profile, the analytical choices associated with each configuration can be displayed alongside the curve. This facilitates the identification of factors that tend to appear in high-performing or low-performing regions and helps reveal structured patterns within the configuration space.

\subsubsection*{Variance Decomposition}

To quantify the relative influence of different analytical decisions on the observed scores, we use variance decomposition based on Analysis of Variance (ANOVA) \citep{stat}. This approach partitions the total variability of a response variable into contributions attributable to specific factors.

In this framework, the response variable corresponds to the score obtained for each evaluated configuration, whereas the explanatory factors represent the analytical choices associated with that configuration. The resulting decomposition estimates how much of the total variability is explained by each factor individually, thereby identifying the main drivers of performance variation.

When the residual component remains large, the analysis can be extended by incorporating interaction terms between selected factors. This makes it possible to assess whether the observed variability is associated with structured dependencies between analytical decisions rather than with isolated main effects alone.

\subsection*{\textbf{Dataset}}

The samples used in this study come from the dataset described in \cite{computacional, nino2025automatic}, where bone marrow flow cytometry data from pediatric B-cell Acute Lymphoblastic Leukemia (B-ALL) were analyzed in the context of relapse prediction. In that work, the authors developed a computational immunophenotyping pipeline based on diagnostic flow cytometry and evaluated its ability to distinguish between patients with and without relapse. The \textit{.fcs} files were acquired and preprocessed following standard clinical flow cytometry procedures, including compensation, transformation and merging of measurements obtained from different tubes \citep{Orfao2002,VanLochem2004,Wood2015}. The data used in the present study are publicly available at \url{https://github.com/Almr95/Relapse-Prediction/tree/main/Selection_A}.

B-ALL is the most common childhood cancer and is characterized by the abnormal proliferation of immature B-cell precursors in the bone marrow. In this context, relapse remains one of the main determinants of poor prognosis, which motivates the distinction between patients with  and without relapse.

Each sample is represented by a matrix
\[
X \in \mathbb{R}^{N \times M},
\]
where each of the \(N\) rows corresponds to an individual cell and each of the \(M\)=8 columns represents an immunophenotypic marker: CD10, CD19, CD20, CD34, CD38, CD45, CD58 and CD66. These specific markers were selected to maintain consistency with the study in which the patient cohort was originally described \citep{computacional}, as they provide the necessary signal for identifying the phenotypic identity of the blasts. These proteins define the $M$-dimensional space where each cell is represented as an observation and whose collective geometry reflects the biological state of the sample at the time of diagnosis. Thus, each sample can be interpreted as a point cloud in \(\mathbb{R}^M\), whose intensity vectors characterize the phenotypic identity of each individual cell.

For this study, we selected patients classified as intermediate risk and treated according to the SEHOP--Pethema 2013 protocol, resulting in a final cohort consisting of diagnostic flow cytometry samples from 76 patients without relaps (NR) and 20 patients with relapse (R). This restriction was introduced to avoid the additional variability that could arise from combining patients managed under different treatment schemes. Within the resulting subset, intermediate-risk patients accounted for approximately 70\% of the relapse cases available in the analyzed dataset. The clinicopathological characteristics of the patients are summarized in Table~\ref{tab:patients}.

\section*{Results}
\label{sec:parameter_selection_results}

This section presents the results obtained from the analysis of the Acute Lymphoblastic Leukemia (ALL) cohort described in Table~\ref{tab:patients}, which is used here as a case study.
 
First, different configurations for the construction of topological representations are compared based on their discriminative capacity between the Non Relapse (NR) and  Relapse (R) groups, evaluated using the Mann-Whitney U (MWU) test. Subsequently, the selected representations are employed as input to examine various supervised classification model configurations and identify those providing the highest predictive performance, measured by the $F_2$ metric.\\

The results are organized into two complementary analytical levels. The first result level, $R_1$, employs a statistical analysis based on the MWU test to select the configuration that best discriminates between patient cohorts for each topological representation. This process identifies the most effective parametric combinations for each representation across homological dimensions $H_0$, $H_1$ and the concatenated $H_0{+}H_1$ vector, which is formed by the combination of the previous two.\\

The second result level, $R_2$, evaluates the performance of classification configurations built from these selected representations. This stage incorporates a global comparison across different representations, homological dimensions and classifiers, while also considering training scenarios with and without oversampling. Together, these analyses make it possible to identify representation settings with strong discriminative power and to determine which modeling combinations yield the best predictive performance.

\(R_1\) and \(R_2\) serve different purposes within the workflow. \(R_1\) acts as a topological screening stage based on statistical separability, whereas \(R_2\) evaluates the predictive utility of the selected representations in a supervised setting.

Although cross-validation is used in both result levels, it serves different purposes in each case. In \(R_1\), it is used to assess the stability of topological configurations in separating the NR and R groups across patient partitions. In \(R_2\), it is used to evaluate classifier performance and generalization on unseen data.

The analysis was conducted sequentially on the same cohort. In \(R_1\), we identified the topological configurations with the highest statistical separability between the NR and R groups. These configurations were then used as input for \(R_2\), where their predictive performance was assessed through supervised classification. Under this design, the results from \(R_2\) should be interpreted as exploratory and comparative within the analyzed cohort, rather than as a fully independent validation of predictive performance.

\begin{figure}[H]
        \centering
        \includegraphics[width=0.99\textwidth]{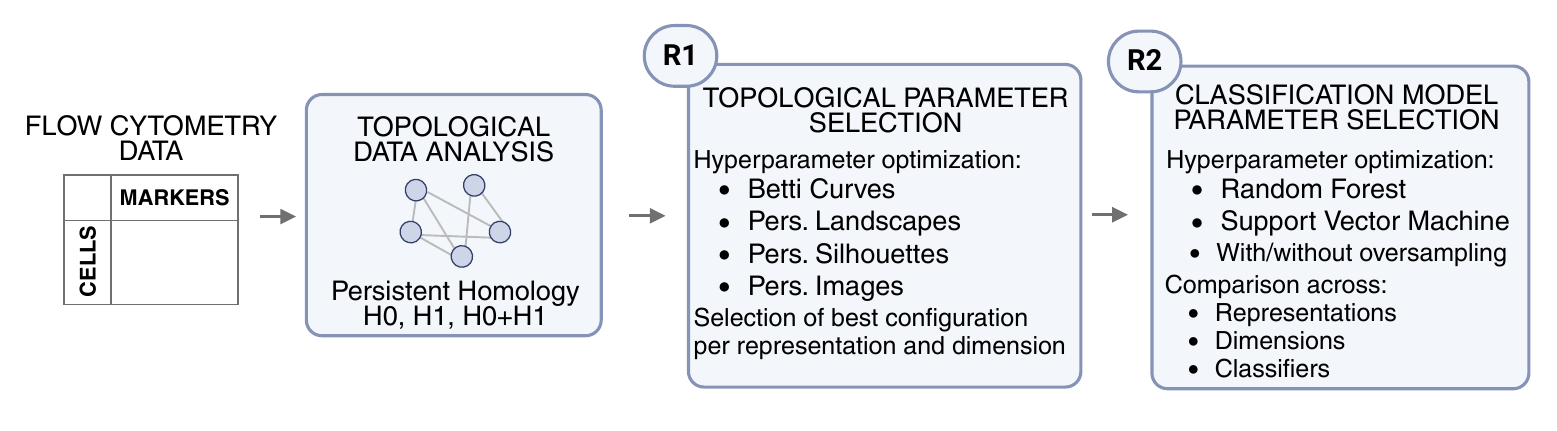}
           \caption{\textbf{Schematic of a two-level analysis framework (\(R_1\) and \(R_2\)).} The workflow integrates the extraction of topological descriptors using persistent homology and systematic optimization. Level \(R_1\) focuses on the selection of representation hyperparameters based on statistical discriminative power, while \(R_2\) evaluates predictive performance using supervised classification models under different modeling configurations.}
            
        \label{fig:results}
    \end{figure}

 \subsection*{$R_1$. Topological Parameter Selection}
\label{subsec:R1_topological_parameter_selection}

\begin{figure}[H]
        \centering
        \includegraphics[width=0.99\textwidth]{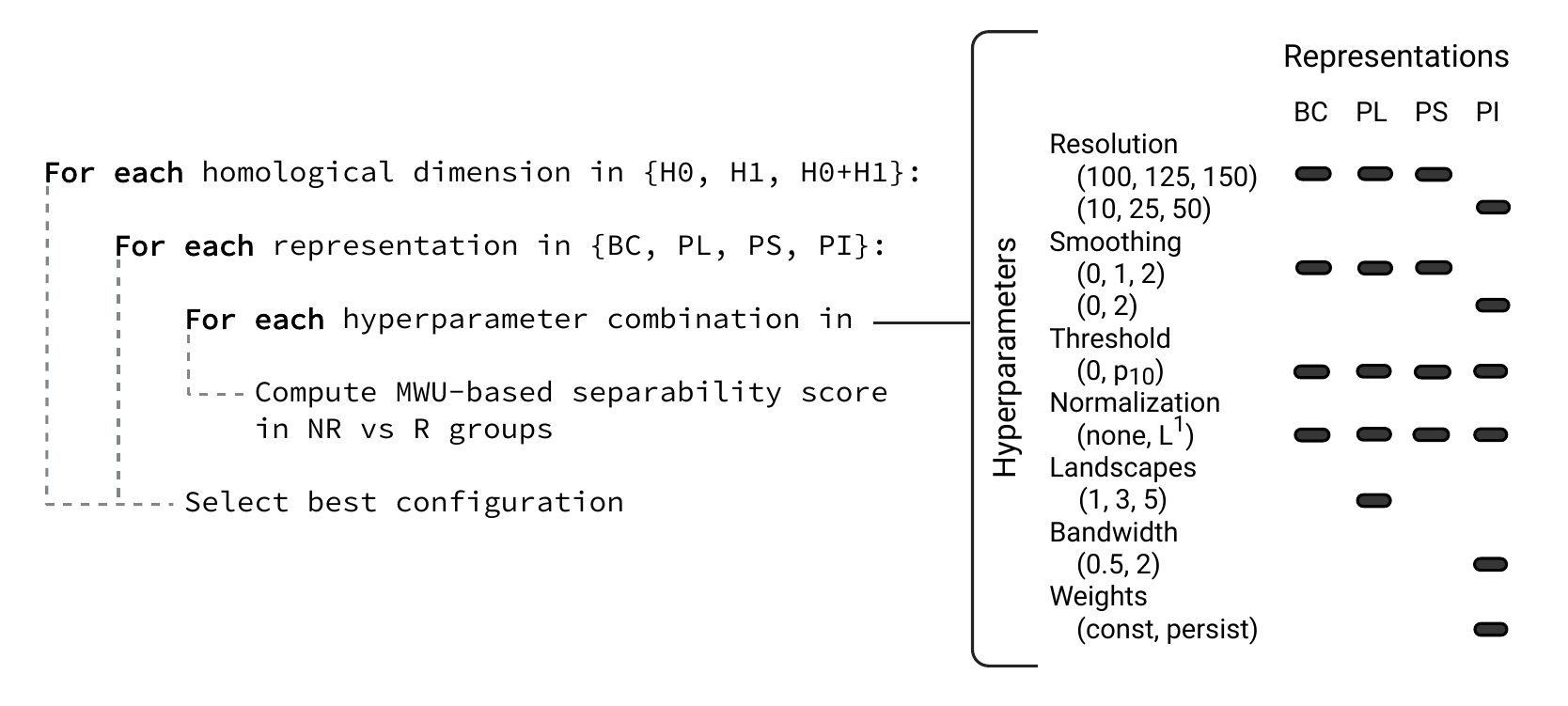}
            \caption{\textbf{Schematic representation of the procedure followed in \(R_1\) for topological parameter selection.} For each homological dimension (\(H_0\), \(H_1\) and \(H_0+H_1\)) and for each topological representation (Betti Curves (BC), Persistence Landscapes (PL), Persistence Silhouettes (PS) and Persistence Images (PI)), all valid hyperparameter combinations are explored. Each configuration is evaluated using an MWU-based separability score computed from the comparison between the NR and R groups. The best configuration is then selected for each representation and homological dimension.}
        
        \label{fig:r1}
    \end{figure}

For each of the four topological representations considered, we evaluated the capacity to discriminate between patients without relapse (NR) and with relapse (R) using the Mann-Whitney U (MWU) test. In each case, we analyzed homological dimension $H_0$, $H_1$ and the concatenated $H_0{+}H_1$ representation, exploring a wide range of hyperparameter combinations to identify configurations that maximize statistical separation (Figure \ref{fig:r1}).

For each evaluated configuration, the MWU-based separability analysis was performed within a stratified cross-validation scheme at the patient level. In each fold, data-dependent quantities, such as the persistence threshold p10, were estimated using only the training partition and then applied to the corresponding test partition. The resulting separability scores were computed on the test folds and summarized across partitions by their median values. Configurations were ranked according to their median score, so that the selected settings correspond to those showing the most consistent separation between the NR and R groups across different data splits.

Specification curves provide an overview of how discriminative power changes across the configurations evaluated for each representation family. They show that the effect of representation choices is not the same across families, and that the role of homological dimension depends on the specific representation. Overall, two main patterns emerge: one in which homological dimension is the dominant factor, and another in which variability is spread across several parameters.

\begin{figure}[H]

    \centering
    \includegraphics[width=1\textwidth]{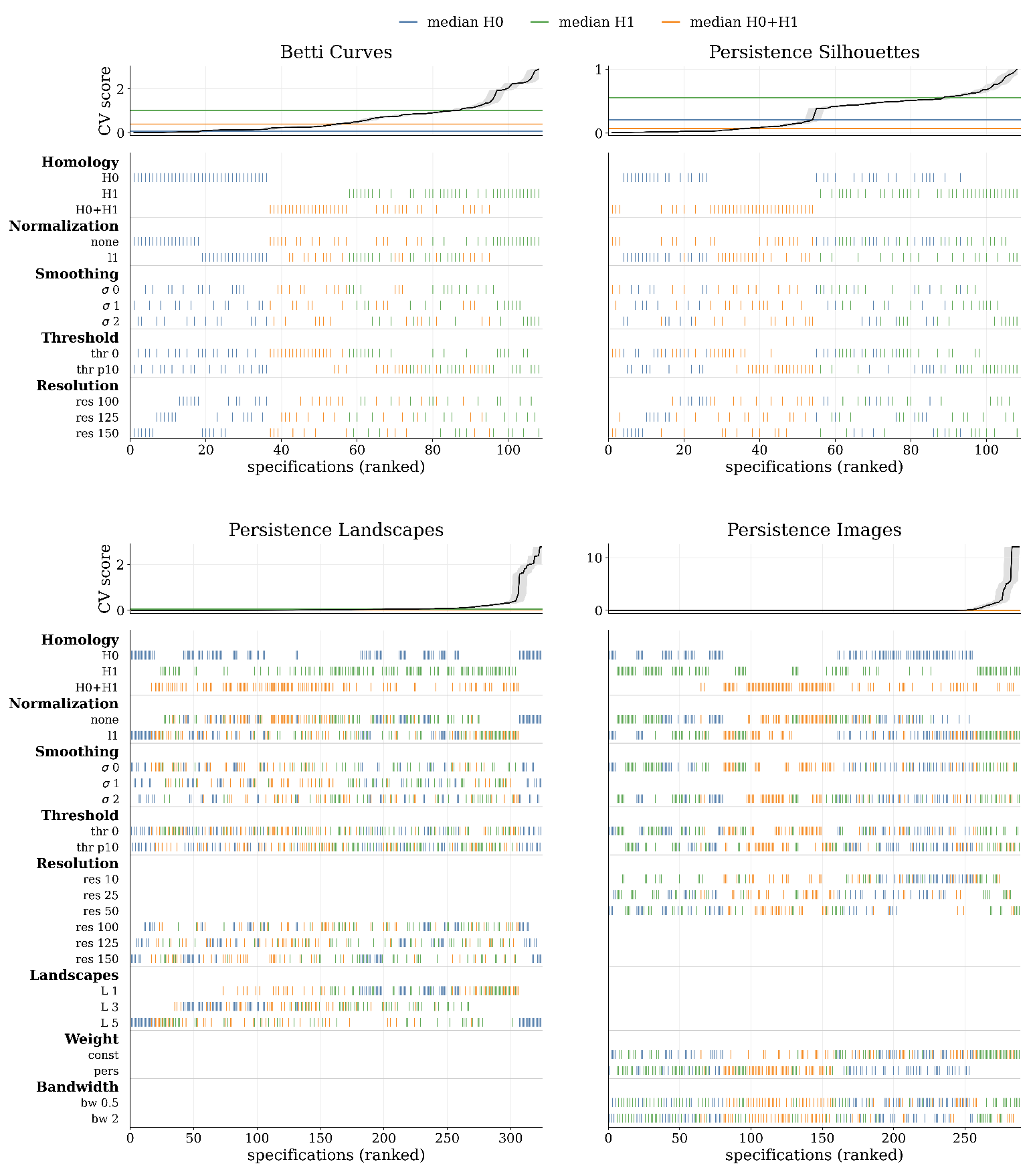}
    \caption{\textbf{Specification curves for the topological representations.}
Specification curves comparing model performance across all evaluated hyperparameter configurations for Betti Curves, Persistence Silhouettes, Persistence Landscapes and Persistence Images. Within each representation, specifications are ranked from lowest to highest median cross-validation score. The black curve shows the performance of each specification and the shaded band indicates the local variability across neighboring specifications. Horizontal colored lines indicate the median performance for each homology setting: $H_0, H_1,$ and $H_0{+}H_1$. The lower panel for each representation shows the hyperparameter choices associated with each ranked specification, with tick colors indicating the corresponding homology dimension. Empty rows indicate hyperparameters that are not applicable to a given representation.}

    \label{fig:spec_curve_mwu}
\end{figure}

As shown in Figure~\ref{fig:spec_curve_mwu}, the four topological representations differ in their sensitivity to the analytical choices considered. For Betti Curves and Persistence Silhouettes, performance changes across configurations follow a relatively structured pattern, suggesting that some representation choices have a clear impact on the discrimination between the NR and R groups. In both cases, homological dimension appears to be one of the main factors driving the observed behavior.

For Betti Curves (BC), the best-performing configuration used \(H_1\), no normalization, a resolution of 150, smoothing with \(\sigma = 2\) and a p10 threshold. For Persistence Silhouettes (PS), the best result was also obtained with \(H_1\), but with \(L^1\) normalization, a resolution of 150, smoothing with \(\sigma = 2\) and a p10 threshold. This consistency suggests that, for these two representations, one-dimensional homological features are particularly informative for distinguishing between the two groups.

Persistence Landscapes (PL) and Persistence Images (PI) behave differently. In PL, the specification curve still suggests a structured pattern, although less clearly dominated by a single factor than in BC and PS. This suggests that discriminative performance in PL depends on a combination of analytical choices rather than on a single dominant parameter. By contrast, PI shows a more dispersed pattern across specifications, with greater sensitivity to the analytical configuration and without a similarly clear dominant factor. PI also spans the widest score range among all representations and reaches the highest separability values in \(R_1\), although this behavior is accompanied by stronger sensitivity to parameterization.

For PL, the best-performing configuration used $H_0$, no normalization, a resolution of 150, smoothing with $\sigma$ = 1, threshold 0 and five landscapes. For PI, the optimal configuration used $H_1$, $L^1$ normalization, a resolution of 50, smoothing with $\sigma$ = 0, threshold 0, a bandwidth of 0.5 and constant weighting. Overall, these results suggest that homological dimension is an important source of variation across all representations, although its relative influence depends on how each representation encodes the underlying topological information.

\begin{figure}[H]
    \centering
    \includegraphics[width=0.99\textwidth]{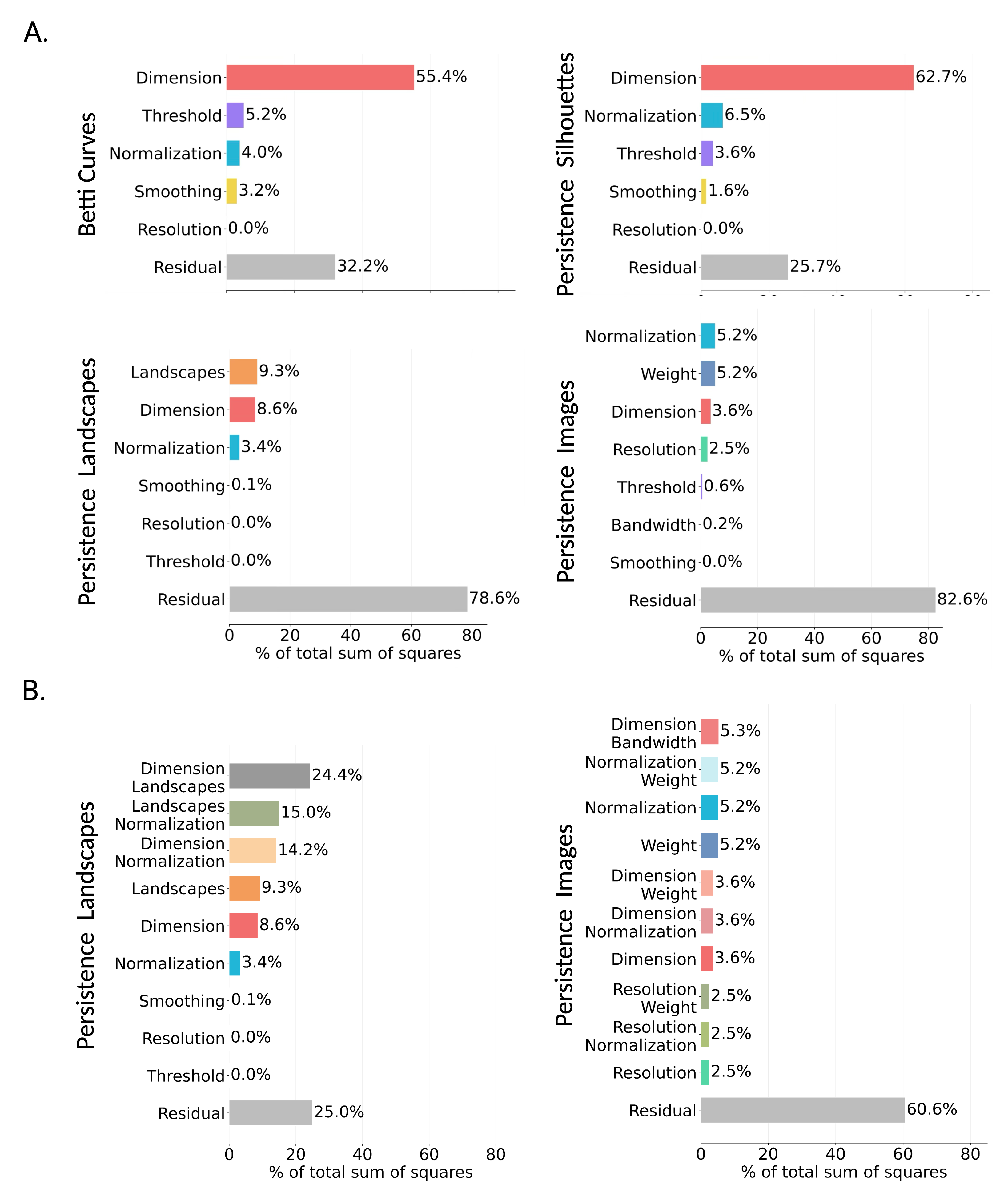}

    \caption{\textbf{Variance decomposition of the separability score across topological representations.} ANOVA-based decomposition of the variability associated with the MWU-derived separability score. \textbf{A.} Main effects, highlighting the dominant role of homological dimension in Betti Curves and Persistence Silhouettes, in contrast to the more distributed variability observed in Persistence Landscapes and Persistence Images. \textbf{B.} Extended analysis for the representations with the largest residual components, incorporating selected second-order interactions. In Persistence Landscapes, these interactions account for a substantial part of the variability, whereas in Persistence Images a large residual term remains, indicating stronger sensitivity to parameterization.}
    \label{fig:mwu_anova}
\end{figure}

Figure~\ref{fig:mwu_anova}A confirms the patterns already observed in the specification curves and helps quantify their relative importance. In BC and PS, homological dimension is the main source of variability in the aggregate score. For BC, dimension explains 55.4\% of the total sum of squares, well above thresholding (5.2\%), normalization (4.0\%), smoothing (3.2\%), and resolution (0.0\%). A similar pattern appears in PS, where dimension accounts for 62.7\% of the total variability, clearly above normalization (6.5\%), thresholding (3.6\%), and smoothing (1.6\%). In both cases, the residual term remains relatively low (32.2\% in BC and 25.7\% in PS), suggesting that a substantial part of the observed variability is captured by the factors included in the model.

PI and PL show a more distributed pattern. In PI, the residual term accounts for 82.6\% of the total sum of squares, which suggests that most of the variability is not explained by the factors explicitly included in the model. Among those factors, normalization and weight are the largest contributors, both at 5.2\%, followed by dimension at 3.6\%. For PL, the largest contributions come from the number of landscapes (9.3\%) and dimension (8.6\%), while the residual term also remains high at 78.6\%. These results are consistent with the specification curves: BC and PS are mainly driven by homological dimension, whereas PL and PI depend on a broader combination of analytical factors.

For the representations with the largest residual components, the analysis was extended by incorporating selected second-order interactions in the same ANOVA framework (see Figure~\ref{fig:mwu_anova}B). In PL, the inclusion of interactions markedly reduces the residual component, from 78.6\% to 25.0\%. This suggests that much of the variability observed in the additive model is associated with structured dependencies between parameters. The most important contributions arise from the dimension-landscapes interaction (24.4\%), followed by landscapes-normalization (15\%) and dimension-normalization (14.2\%). In PI, the reduction is more moderate, from 82.6\% to 60.6\%. The most relevant interactions are dimension-bandwidth (5.3\%) and normalization-weight (5.2\%). In this case, the interaction terms capture only part of the observed heterogeneity, which is consistent with a stronger and less stable dependence on parameterization.

Taken together, the variance decomposition shown in Figure~\ref{fig:mwu_anova} refines the interpretation of \(R_1\). BC and PS show a more stable structure, with homological dimension as the main explanatory factor. PL also exhibits a structured pattern, although its behavior depends more clearly on combinations of analytical choices, particularly interactions involving dimension, normalization, and the number of landscapes. PI, by contrast, reaches the highest score range in \(R_1\) while also remaining the most sensitive to hyperparameter selection, as a substantial part of its variability remains unexplained even after incorporating second-order interactions. In this way, \(R_1\) helps identify configurations with strong separation between the NR and R cohorts, while also distinguishing between robust representations and those that depend more strongly on parameterization.

The statistical screening performed in Result Level \(R_1\) is complemented by a detailed visualization of the performance grids in the Supplementary data. Specifically, Figures S1–S4 document the MWU-based score distributions across the entire hyperparameter space for each topological family, while Figure S5 provides an extended variance decomposition analysis stratified by homological dimension.

Accordingly, for the subsequent \(R_2\) analyses, we retained the best-performing configuration within each homological setting and representation. This ensured that comparisons across \(H_0\), \(H_1\) and \(H_0+H_1\) reflected the strongest available specification for each representation, rather than being driven by a single global optimum.

\subsection*{$R_2$. Classification parameter selection}
\label{subsec:R2_classification_parameter_selection}

This subsection evaluates the performance of supervised classification models constructed from the topological representations selected in \(R_1\). For each representation and homological dimension (\(H_0\), \(H_1\) and the concatenated \(H_0+H_1\) vector), we considered Random Forest (RF) and Support Vector Machine (SVM) classifiers. We explored different hyperparameter combinations under two training scenarios, with and without oversampling of the minority class. This analysis helps assess how topological representation, homological dimension, classifier family and balancing strategy affect predictive performance, using \(F_2\) as the main metric.

For each homological dimension, oversampling setting, and classifier hyperparameter combination, evaluation was performed within a stratified cross-validation scheme at the patient level. In each fold, preprocessing steps were estimated using only the training partition: the StandardScaler was fitted on the training data, oversampling was applied only to the training partition when required, and the classifier was then trained on the training fold. Performance was then evaluated on the corresponding test partition.

For each configuration, AUC, accuracy, and \(F_2\) were averaged across folds. In addition, out-of-fold predictions were obtained to construct an aggregated confusion matrix for each evaluated model. This procedure ensures that scaling and balancing are recalculated independently within each fold, which avoids information leakage from the test partitions into the training process.

\begin{figure}[H]
        \centering
        \includegraphics[width=0.99\textwidth]{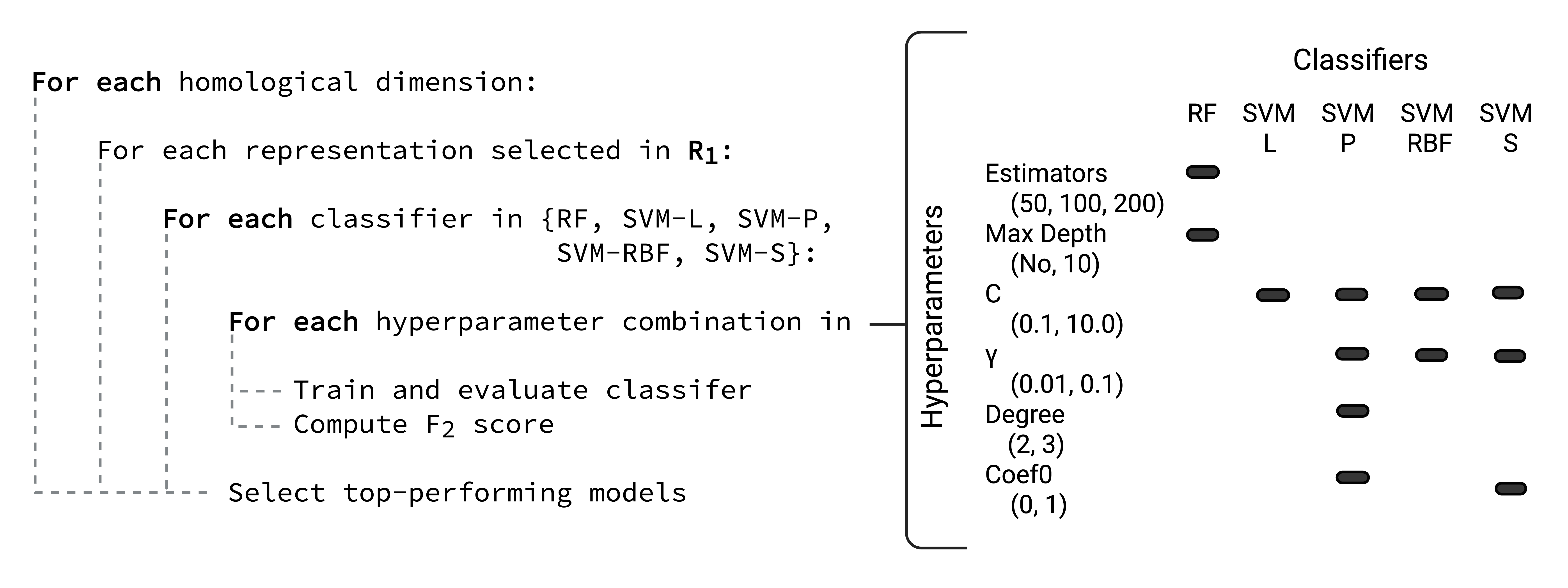}
        \caption{\textbf{Schematic representation of the procedure followed in \(R_2\) for classification parameter selection and model assessment.} For each homological dimension and for each topological representation (Random Forest (RF), Support Vector Machine Linear Kernel (SVM L), Polynomial Kernel (SVM P), Radial Basis Function Kernel (SVM RBF) and Sigmoid Kernel (SVM S)) selected in \(R_1\), all valid classifier and hyperparameter combinations are explored. Each configuration is trained and evaluated using the \(F_2\) score, and the top-performing models are retained. These models are then further assessed using AUC and the corresponding confusion matrices to identify the most suitable configurations for classifying the NR and R groups.}
        \label{fig:r2}
    \end{figure}

\begin{figure}[H]
\centering
    \includegraphics[width=1\textwidth]{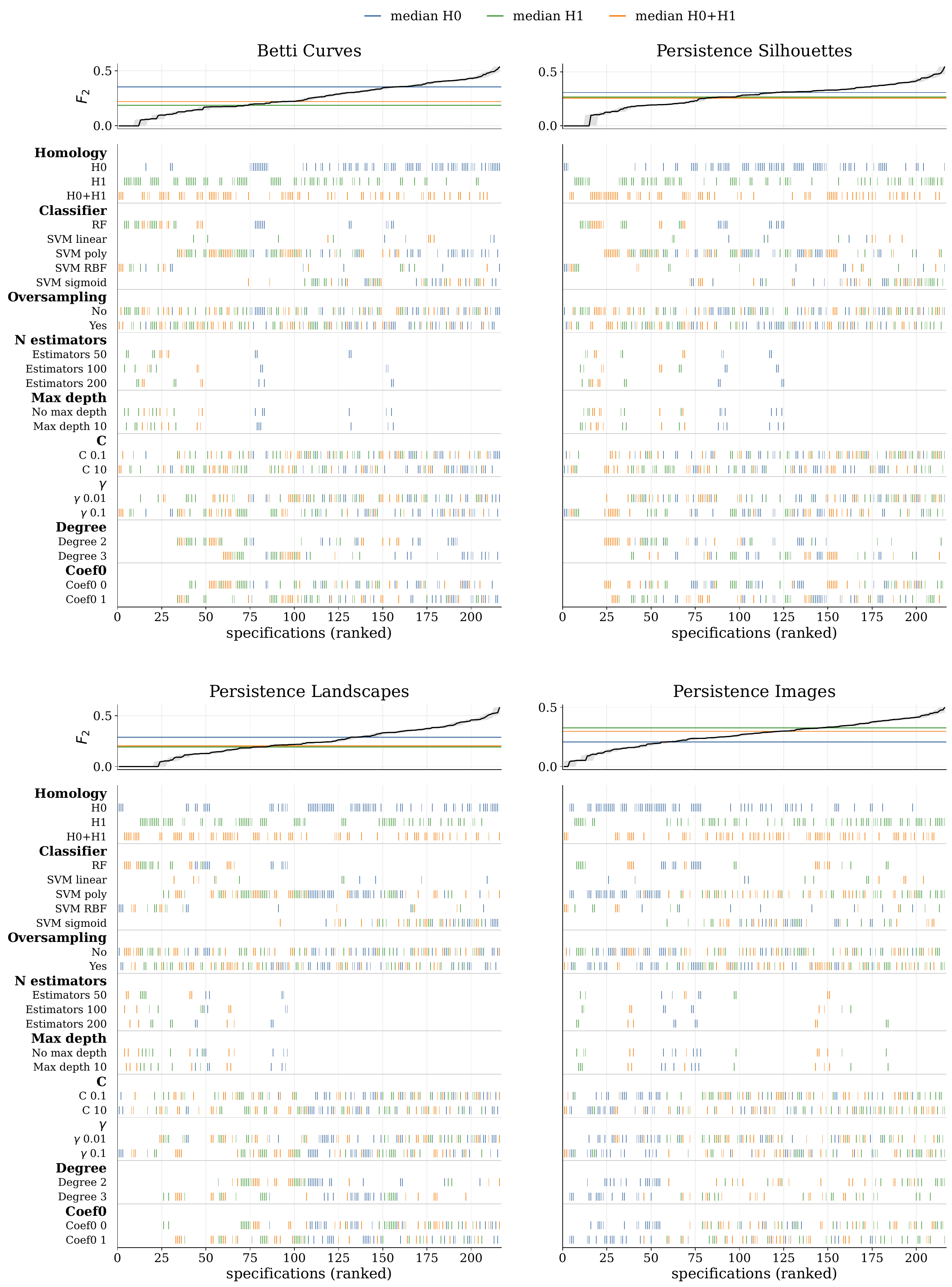}
    \caption{\textbf{Specification curves for the \(R_2\) classification analyses across the four topological representations.} Curves show the \(F_2\) score across all evaluated classification configurations for Betti Curves, Persistence Silhouettes, Persistence Landscapes and Persistence Images, ranked from lowest to highest performance within each representation. The black curve denotes specification-level performance, the shaded band indicates local variability and horizontal colored lines show the median \(F_2\) score for \(H_0\), \(H_1\) and \(H_0{+}H_1\). Lower panels display the corresponding classifier, oversampling strategy and model-specific hyperparameters; tick colors indicate homological dimension and empty rows denote non-applicable parameters.}
    \label{fig:spec_curve_f2}
\end{figure}

As shown in Figure~\ref{fig:spec_curve_f2}, the \(R_2\) specification curves provide a compact view of how supervised classification performance varies across the full set of evaluated configurations. In contrast to \(R_1\), which focused on statistical separation between groups, \(R_2\) examines how these topological summaries perform when combined with different classifiers, oversampling strategies and model-specific hyperparameters. Overall, the results show that classification performance is determined by the interplay between the topological representation, the learning algorithm, and the specific hyperparameter settings.

For BC, the highest-performing configurations are mainly in \(H_0\), whose median \(F_2\) score is higher than those of \(H_1\) and \(H_0{+}H_1\). The best BC configuration used \(H_0\), a linear SVM, no oversampling and \(C = 0.1\). This suggests that, for this representation, relatively simple decision boundaries applied to connected-component information provide the best classification results.

PS show a more balanced pattern across homological dimensions, although \(H_0\) remains among the most competitive options. The best-performing PS configuration used \(H_0\), an RBF SVM, no oversampling, \(C = 0.1\) and \(\gamma = 0.01\). Compared with BC, this suggests a slightly more nonlinear classification structure, while still relying primarily on \(H_0\)-based topological information.

PL, performance is more distributed across the configuration space, with \(H_0{+}H_1\) yielding the best overall result. The optimal PL configuration used \(H_0{+}H_1\), a polynomial SVM, no oversampling, \(C = 0.1\), \(\gamma = 0.01\), degree 2 and \(\mathrm{coef0} = 0\). This suggests that, for landscapes, combining information from both homological dimensions can be beneficial when paired with a classifier able to capture nonlinear interactions.

PI display a distinct pattern, with the best-performing configuration based on \(H_1\). The optimal PI configuration used \(H_1\), a polynomial SVM, oversampling, \(C = 10\), \(\gamma = 0.1\), degree 2 and $coef0 = 1$. This suggests that PI is more sensitive to classifier complexity and parameterization, which is consistent with the higher sensitivity to specification choices observed for this representation. Overall, the \(R_2\) results indicate that the most informative homological dimension depends on the representation and that the best classification performance comes from specific combinations of topological summary, classifier family and hyperparameter setting.

\begin{figure}[H]
    \centering
    \includegraphics[width=0.99\textwidth]{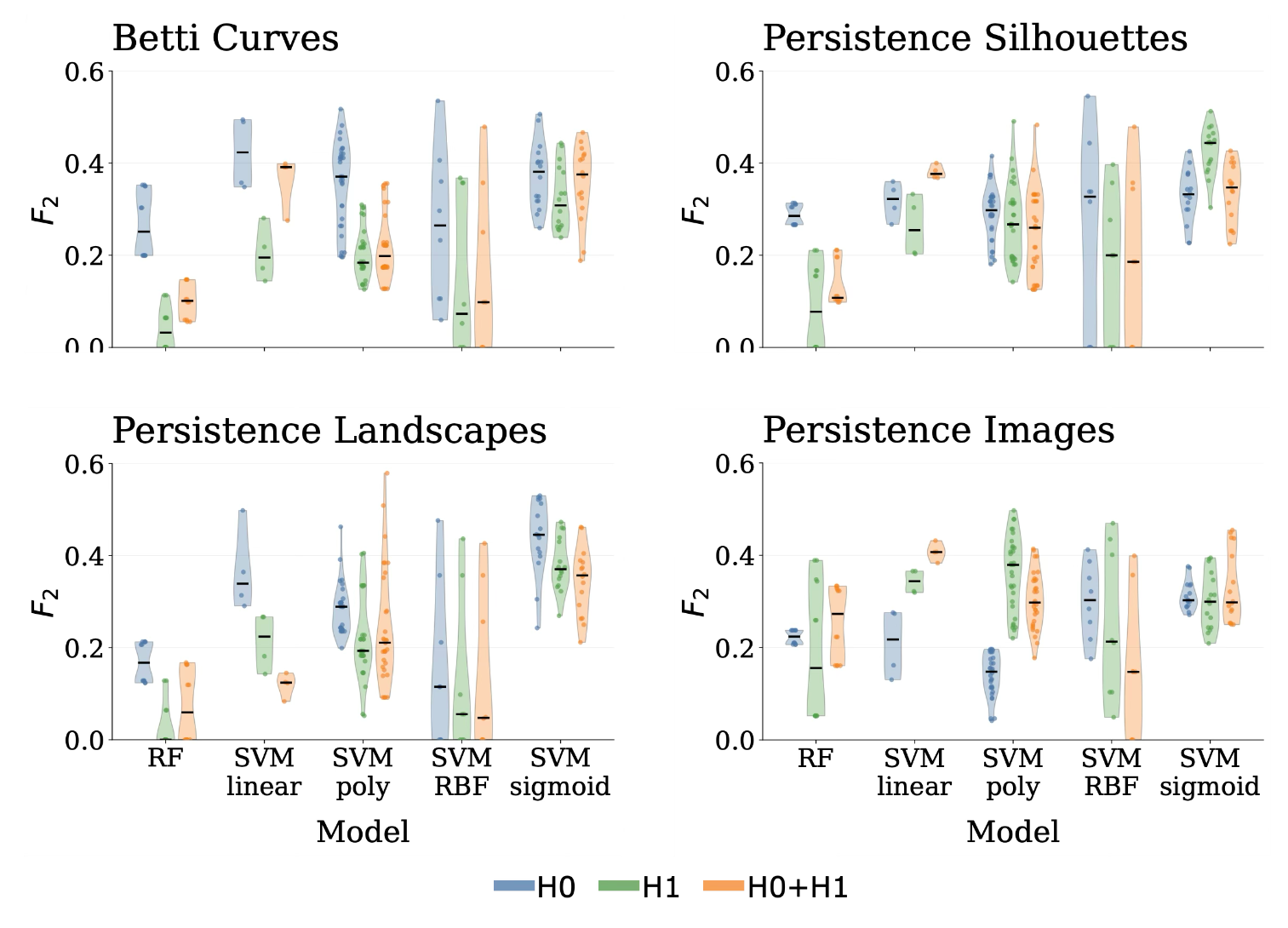}

        \caption{\textbf{Marginal distributions of $F_2$ by classifier and homological dimension.} Comparative analysis of supervised performance illustrating the interaction between classifiers and topological representations.}

    \label{fig:f2_marginal_violin}
\end{figure}

The marginal distributions of \(F_2\) scores shown in Figure~\ref{fig:f2_marginal_violin} provide a more detailed view of the differences across representations, dimensions and classifier families. In BC, the highest scores are mainly concentrated in linear and sigmoidal SVM models, especially for the \(H_0\) and \(H_0{+}H_1\), whereas the RBF kernel SVM shows significantly greater variability. In PI, the \(H_1\) remains competitive across several models, particularly with polynomial and sigmoidal SVM kernels. 

For PL, the sigmoidal SVM again reaches some of the highest scores, specifically in the \(H_0\) dimension. PS, by contrast, shows a more balanced pattern, with relevant contributions from \(H_0\), \(H_1\) and \(H_0{+}H_1\) depending on the classifier used. Overall, these distributions suggest that there is no universally optimal representation, since the usefulness of each homological dimension depends on its interaction with the classifier family.

\begin{figure}[H]
    \centering
    \includegraphics[width=0.99\textwidth]{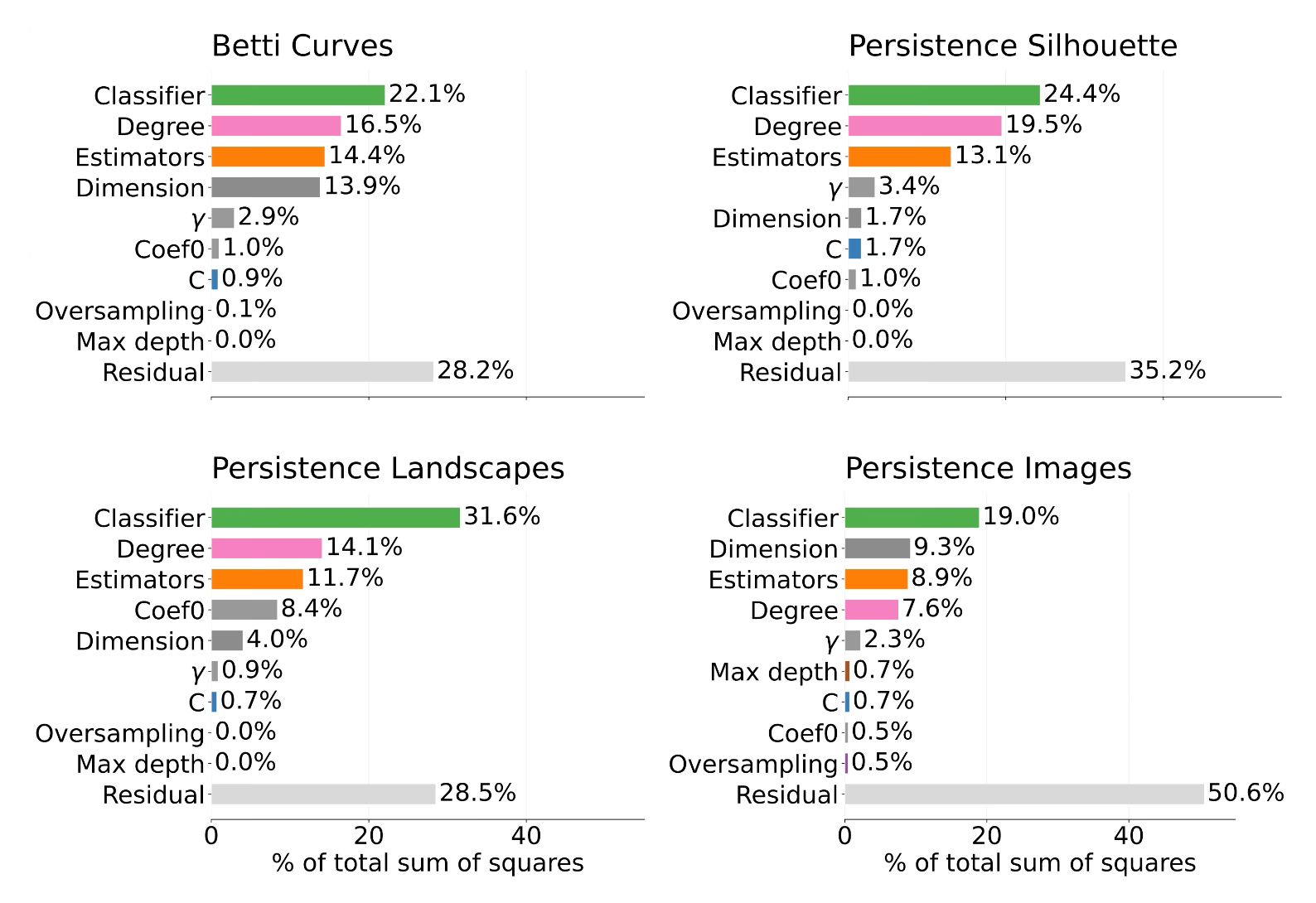}

      \caption{\textbf{Variance decomposition based on ANOVA for the $F_2$ score across the four families of topological representations.} Analysis of Betti Curves, Persistence Silhouettes, Persistence Landscapes and Persistence Images. In Betti Curves and Persistence Silhouettes, the classifier constitutes the primary explanatory factor for the observed variability.}

    \label{fig:class_anova}
\end{figure}

The variance decomposition presented in Figure~\ref{fig:class_anova} quantifies the patterns previously observed in the specification curves and the marginal distributions of \(F_2\) scores. In BC, the classifier family is the main explanatory factor, accounting for 22.1\% of the total variability, followed by the polynomial degree (16.5\%), the number of estimators  (14.4\%) and the homological dimension (13.9\%). In PI, the classifier also remains the dominant factor at 19.0\%, although this representation shows a much higher residual fraction of 50.6\%. 

For PL, the classifier explains 31.6\% of the variability, followed by degree (14.1\%), the number of estimators (11.7\%) and the $coef0$ (8.4\%). A similar pattern appears in PS, where the classifier accounts for 24.4\% of the variability, followed by the degree (19.5\%) and the number of  estimators (13.1\%). Overall, these results show that classifier family has a strong influence on supervised performance. Its effect, however, depends on how it interacts with other modeling factors in each representation.

\begin{figure}[H]
    \centering
    \includegraphics[width=0.65\textwidth]{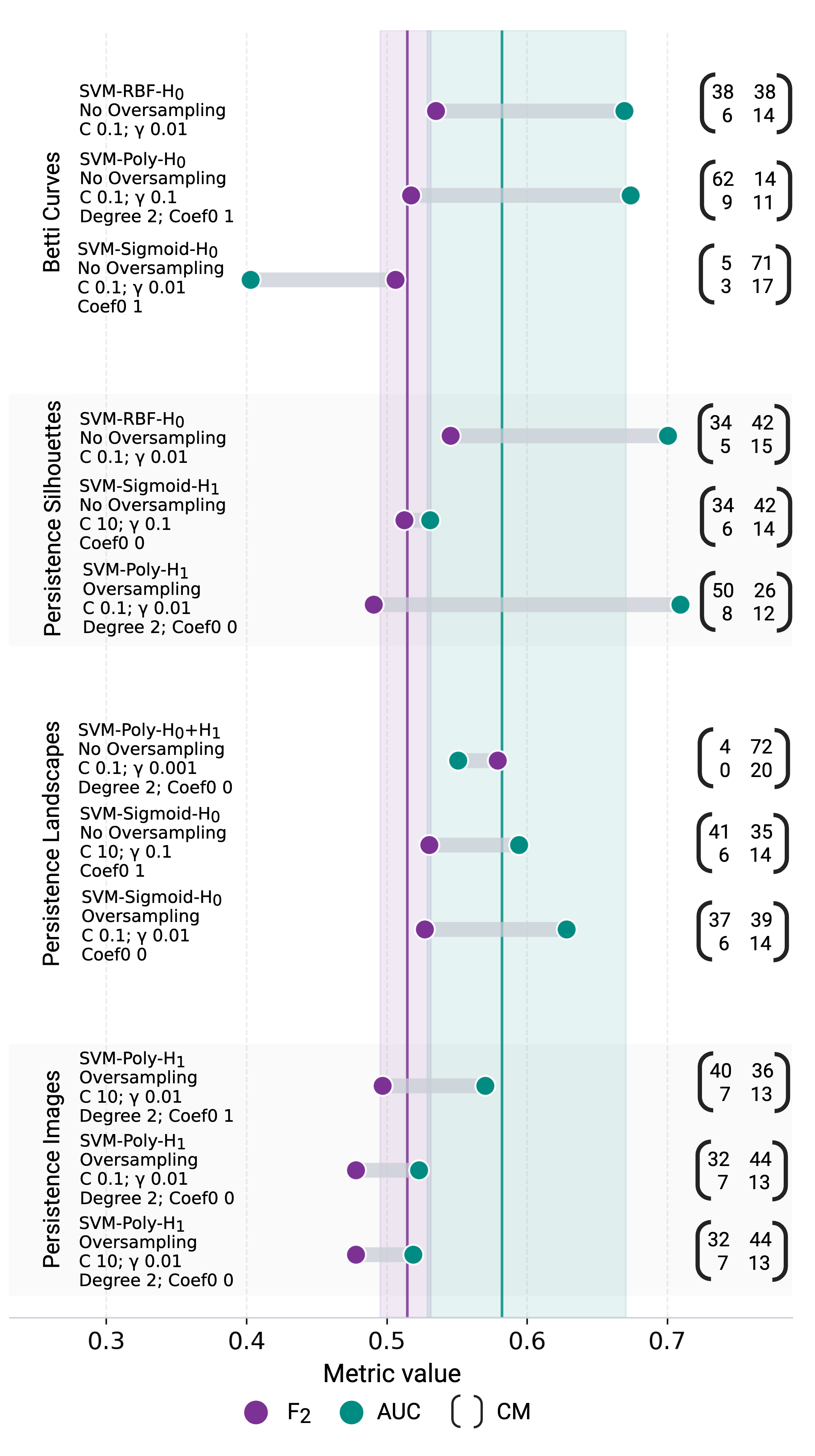}

        \caption{\textbf{Top three classifier configurations ranked by \(F_2\) within each topological representation.} For each representation, the three models with the highest \(F_2\) scores are shown. Each horizontal segment connects the \(F_2\) score and AUC of the same model, enabling a joint comparison of performance for the positive class and overall discriminative ability. Vertical lines indicate the median value of each metric across the selected models, and the shaded bands represent the corresponding percentile range. Purple markers and shading denote \(F_2\), whereas green markers and shading denote AUC.}

    \label{fig:f2_auc}
\end{figure}

Figure~\ref{fig:f2_auc} jointly compares \(F_2\), AUC and the confusion matrices for the three best models according to \(F_2\) in each representation. This analysis summarizes the aggregate performance of each configuration and also shows how errors are distributed between the NR and R classes.

The effect of oversampling is not uniform across models. In Random Forest, oversampling the minority class tends to increase the \(F_2\), suggesting an improvement in the detection of relapse cases, while accuracy remains at reasonable levels. This pattern appears across several representations, particularly in Persistence Silhouettes and Persistence Images. In these cases, the increase in correctly identified relapses is not accompanied by a drastic loss in the correct classification of patients without relapse.

In SVM, the effect of oversampling is more variable. Some combinations improve both \(F_2\) and AUC, whereas others lead to less balanced models. In those cases, higher sensitivity is accompanied by substantial losses in accuracy or by a stronger bias toward relapse prediction. This suggests that the effect of the balancing strategy depends on its interaction with both the classifier family and the underlying topological representation.

More generally, AUC values tend to be higher than \(F_2\) scores. This suggests that several models have global discriminative capacity that does not always translate into balanced performance for the positive class. The confusion matrices also show that models with similar \(F_2\) or AUC scores can differ significantly in the balance between false positives and false negatives.

In this context, the second model in the Betti Curves panel of Figure~\ref{fig:f2_auc} appears to provide one of the most favorable trade-offs. It is based on a polynomial-kernel SVM without oversampling in the \(H_0\) dimension. Although it does not achieve the maximum value in any aggregate metric, it reaches an \(F_2\) of 0.52 and an AUC of 0.67, while correctly classifying 62 NR cases and 11 R cases. This suggests a reasonable balance between the detection of relapse cases and the preservation of correct classification among patients without relapse.

In contrast, some models with similar or higher AUC values show less balanced error distributions. This is reflected either in a marked increase in false positives or in insufficient detection of the R class. Figure~\ref{fig:f2_auc} therefore shows that the joint evaluation of \(F_2\), AUC and the confusion matrix is necessary to identify configurations with a truly useful balance between relapse sensitivity and overall classification stability.

The difference observed for PI between the two levels of analysis is particularly striking. In \(R_1\), this representation reaches some of the highest statistical separability scores. However, this does not translate directly into predictive performance in \(R_2\), where PI is outperformed by Betti Curves and Persistence Silhouettes. Sensitivity analysis helps explain this result, since PI shows a very high residual term (82.6\%) in \(R_1\). This suggests that its discriminative signal is highly sensitive to specific and unstable parametric settings. In the supervised context of \(R_2\), this variability acts as noise for the classifiers and limits generalization, especially when compared with more stable descriptors such as BC and PS, where variability is better captured by explicit factors like the homological dimension.

A comprehensive breakdown of the predictive performance for all evaluated classification models is provided in Tables \ref{tab:bc_rf_svm}–\ref{tab:rf_svm_Images} of the Supplementary data. These tables summarize the AUC, accuracy, and \(F_2\) metrics for every configuration tested in Result Level \(R_2\), with the corresponding variance analysis for the supervised stage presented in Figure \ref{fig:anova_dims}.

Overall, the results from \(R_2\) show that supervised performance depends on the interaction between topological representation, homological dimension, classifier family, hyperparameters and balancing strategy. Among these factors, classifier family is one of the most influential across all representation families, although its role changes depending on the representation and the way it interacts with the other modeling choices.

\section*{Discussion and Conclusions}

In this work, we presented a two-level framework for the systematic evaluation of topological representations in high-dimensional single-cell classification. By separating statistical topological screening (\(R_1\)) from supervised classification analysis (\(R_2\)), the proposed strategy makes it possible to compare representations, homological dimensions, model families, and parameter choices within a single coherent workflow. More importantly, it allows the joint assessment of discriminative performance and robustness to analytical choices. Because both \(R_1\) and \(R_2\) were conducted on the same cohort, the predictive results should be interpreted as exploratory and comparative rather than as externally validated performance estimates.

Our results show that topological representations differ not only in their apparent discriminative capacity but also in their robustness to analytical choices. In particular, Persistence Images achieved high statistical separation in \(R_1\), but their behavior in \(R_2\) was more dependent on specific configurations. By contrast, Betti Curves and Persistence Silhouettes, especially in the \(H_0\) dimension, showed a more stable pattern across analyses. These findings suggest that topological information related to connectivity may be particularly relevant for relapse stratification in this cohort.

This difference becomes especially clear in the case of Persistence Images (PI). While PI achieved high statistical scores in the initial phase, its predictive performance during supervised classification was limited by a residual component of 82.6\% in its discriminative signal. In the context of precision medicine, this sensitivity acts as noise and suggests that the apparent usefulness of some topological tools may be compromised by unstable parametric settings. This makes the identification of robust descriptors such as Betti Curves (BC) and Persistence Silhouettes (PS) methodologically important if topological biomarkers are to generalize to unseen patient data.

The results also show that evaluating TDA in this setting requires a multi-metric approach. Relying on a single score such as \(F_2\) can lead to biased models, as illustrated by the case of Persistence Landscapes with SVM-Poly, which failed to identify most non-relapse cases. More broadly, the proposed framework makes it possible to distinguish between representations that perform well only under very specific settings and those that remain more stable across a broader analytical space.

Several limitations should nevertheless be acknowledged. First, the cohort remains relatively small, which limits the statistical strength of the analysis and the stability of some observed patterns. Comparisons with previous studies should also be made with caution. For example, \cite{bib1} reported stronger predictive performance in pediatric acute lymphoblastic leukemia, but the study design is not directly comparable to ours. Our cohort includes only intermediate-risk patients treated according to the SEHOP--Pethema 2013 protocol, which likely defines a more homogeneous and therefore more challenging classification setting. In addition, some representations associated with high performance in earlier studies, particularly Persistence Images, were found here to be more sensitive to parameter choices and less stable across analytical decisions. This suggests that differences between studies may reflect not only cohort composition, but also the robustness of the topological descriptors themselves.

The present study also has methodological limits that point to natural directions for future work. The analysis was restricted to \(H_0\) and \(H_1\), whereas higher homological dimensions may contain additional structural information, albeit at a substantially higher computational cost. Likewise, the framework currently relies on Vietoris--Rips filtrations and the Euclidean metric, and it would be worth exploring alternative constructions, such as Alpha complexes, as well as other distance functions that may better capture the geometry of these data. Larger cohorts, broader modeling settings, and external validation strategies will also be necessary to further assess the robustness and generalizability of the proposed approach.

Beyond absolute predictive performance, the main contribution of this work is methodological. The proposed framework not only identifies high-performing configurations, but also reveals how sensitive different topological representations are to analytical choices. This makes it possible to distinguish between representations that perform well only under very specific settings and those that remain more stable across a broader configuration space. In this sense, the study contributes to a more systematic strategy for applying TDA in complex biomedical data. The aim of the work was not to provide an exhaustive benchmark against all possible non-topological approaches, but to establish a systematic framework for comparing topological representations under a common analytical scheme. The emphasis is therefore placed on relative robustness, sensitivity to parameterization, and consistency across modeling decisions, rather than on maximizing a single predictive metric in isolation.

A further strength of this study is methodological transparency. The full implementation is openly available at \url{https://github.com/RPiconGonzalez/TDA-LLA}, which supports reproducibility, external validation, and reuse in other biomedical applications.

Although the present results are specific to this pediatric acute lymphoblastic leukemia cohort, the proposed framework is generalizable and may be applied to other high-dimensional biomedical datasets in which structural heterogeneity is expected to play an important role. Future work should focus on validating these findings in independent and multicenter cohorts. It should also explore the biological interpretation of the stable topological patterns identified here.

\section*{Acknowledgements}

This work was supported by project PID2022-140451OA-I00, funded by the Ministerio de Ciencia e Innovación and the Agencia Estatal de Investigación (DOI: 10.13039/501100011033). The Asociación Pablo Ugarte (APU, Spain), Junta de Andalucía (Spain), group FQM-201, and the INiBICA research group CO29.

\appendix

\renewcommand{\thetable}{S\arabic{table}}
\setcounter{table}{0}

\renewcommand{\thefigure}{S\arabic{figure}}
\setcounter{figure}{0}

\section*{Supplementary data}
\label{supplementary}

\begin{table}[H]
     \resizebox{\columnwidth}{!}{%
		\begin{tabular}
  {@{\extracolsep{\fill}}lcccc}
\hline
			&  \begin{tabular}[c]{@{}c@{}}Dataset 1 (HVR)\\ (N=5)\end{tabular} & \begin{tabular}[c]{@{}c@{}}Dataset 2 (HVA) \\ (N=14)\end{tabular} & \begin{tabular}[c]{@{}c@{}}Dataset 3 (HNJ)\\ (N=77)\end{tabular} & \begin{tabular}[c]{@{}c@{}}Total\\ (N=96)\end{tabular} \\
			\hline
			Sex - no. (\%)  &  &  &   &   \\
			\qquad Male            &   1 (20.0) &       8 (57.1) &       36 (46.8) &      45 (46.9) \\
			\qquad Female          &   4 (80.0) &       6 (42.9) &       41 (53.2) &      51 (53.1) \\
			\hline
			Age at diagnosis - yr    &  &  &   &  \\
			\qquad Median          & 3 &              7&              4 &              4\\
			\qquad Range          & 2 - 12 &      2 - 15 &      1 - 15 &      1 - 15\\
			\hline
			Long term status -no. (\%)  &  &  & & \\
			\qquad Relapse        & 1 (20.0) &         2 (14.3) &       17 (22.1) &       20 (20.8)\\
			\qquad No relapse      & 4 (80.0) &       12 (85.7) &       60 (77.9) &      76 (79.2)\\
			\hline
			Immunophenotype - no. (\%)  &  &   &   &   \\
			\qquad Common        &   4 (80.0) &       13 (92.9) &      71 (92.2) &      88 (91.7)\\
			\qquad Pre-B            & 1 (20.0) &       1 (7.1) &         4 (5.2) &       6 (6.2)\\
			\qquad Pro-B             &  0 (0.0) &         0 (0.0) &         2 (2.6) &         2 (2.1)\\

			\hline
			BM blasts at diagnosis - \%  &  &  &   &   \\
			\qquad Median          &  79.1 &             77.1 &             86.0 &             85.0 \\
			\qquad Range           & 69.7 - 89.9 &     46.6 - 94.0 &      30.0 - 99.0 &      30.0 - 99.0 \\   \hline
			
			Leukocytes - cell/nL  &  &  &   &   \\
			\qquad Median          &    5.33 &            7.56 &           10.91 &            9.38 \\
			\qquad Range           & 3.58 - 124.3 &  0.54 - 50.0 &  0.21 - 86780.0 &  0.21 - 86780.0 \\   \hline
			
			CNS - no. (\%) &  &  &   &   \\
			\qquad Yes  			& 0 (0.0) &         1 (7.1) &         6 (7.8) &        7 (7.3) \\
			\qquad No & 5 (100.0) &         13 (92.9) &         71 (92.2) &        89 (92.7) \\
			\hline
			Karyotype - no. (\%) &  &  &   &   \\
			\qquad Hyperdiploid ($>$50) & 0 (0.0) &         2 (14.3) &       10 (13.0) &       12 (12.5)\\
			\qquad Normal (40-50) & 3 (60.0) &       3 (21.4) &       44 (57.1) &       50 (52.1)\\
			\qquad Hypodiploid ($<$40) &   0 (0) &        0 (0) &       0 (0) &       0 (0)\\
			\qquad	 No metaphases & 2 (40.0) &        0 (0.0) &       22 (28.6) &       24 (25.0) \\
               \qquad	 No data &  0 (0.0) &        9 (64.3) &       1 (1.3) &       10 (10.4) \\
			\hline
			Chromosomic alterations - no. (\%) &  &  &   &   \\
			\qquad ETV6/RUNX1 t(12;21) &   2 (40.0) &       2 (14.3) &       21 (27.3) &       25 (26.0)\\
			\qquad TCF3/PBX1 t(1;19) & 0 (0.0) &         0 (0.0) &         4 (5.2) &         4 (4.2) \\
			\qquad MLL rearrangement & 0 (0.0) &         0 (0.0) &          1 (1.3) &         1 (1.0)\\
			\qquad BCR/ABL1 t(9;22) &  0 (0)&        0 (0.0) &       0 (0.0)  &        0 (0.0)\\
               \qquad No data &  3 (60.0)&         12 (85.7) &         51 (66.2) &         66 (68.8)\\
			\hline
		\end{tabular}
        }
		\vspace{0.5cm}		
		\caption{Summary of the clinicopathologic characteristics of the cohort. HVR = Virgen del Rocío Hospital, HVA = Virgen de la Arrixaca Hospital, HNJ = Niño Jesus Hopital. Table adapted from \cite{computacional}.}
        \label{tab:patients}
	\end{table}

\subsection*{\(R_1\)}
\label{supp-r1}

\begin{table}[H]
    \centering
    \small
    \resizebox{\columnwidth}{!}{%
    \begin{tabular}{||l||llll||}
    \hline\hline
        \textbf{Hyperparameter} & \textbf{Betti Curves} & \textbf{P. Landscapes} & \textbf{P. Silhouettes} & \textbf{P. Images} \\ \hline \hline
        \textbf{Resolution} & (100, 125, 150) & (100, 125, 150) & (100, 125, 150) & (10, 25, 50) \\ 
        \textbf{Smoothing} & (0, 1, 2) & (0, 1, 2) & (0, 1, 2) & (0, 2) \\ 
        \textbf{Threshold} & (0, p10) & (0, p10) & (0, p10) & (0, p10) \\ 
        \textbf{Normalization} & (None, $L^1$ ) & (None, $L^1$ ) & (None, $L^1$ ) & (None, $L^1$ ) \\ 
        \textbf{Landscapes} & ~ & (1, 3, 5) & ~ & ~ \\ 
        \textbf{Bandwidths} & ~ & ~ & ~ & (0.5, 2) \\ 
        \textbf{Weights} & ~ & ~ & ~ & (const, persist) \\ 
        \hline\hline
    \end{tabular}%
    }
    \caption{Grid of hyperparameters explored for each topological representation. For each representation, different combinations of parameters were evaluated in order to identify the configurations that maximize the statistical separation between the NR and R groups.}
    \label{grid}
        
\medskip
\footnotesize
\textit{Note:} $p10$ indicates the 10th percentile of persistence, used as a threshold to filter out bars with low persistence. 
$L^1$ denotes $L^1$ norm normalization. 
\textit{const} corresponds to constant weights for all points in the persistence plot, whereas  
\textit{persist} indicates weights proportional to the persistence of each point.

\end{table}

\renewcommand{\arraystretch}{1.4}  

\begin{table}[H]
    \centering
    \resizebox{\columnwidth}{!}{%
    \begin{tabular}{||l||ccc|ccc|ccc|ccc||}
    \hline\hline
    \multirow{2}{*}{\textbf{Hyperparameter}} 
      & \multicolumn{3}{c|}{\textbf{Betti Curves}} 
      & \multicolumn{3}{c|}{\textbf{P. Landscapes}} 
      & \multicolumn{3}{c|}{\textbf{P. Silhouettes}} 
      & \multicolumn{3}{c||}{\textbf{P. Images}} \\
    \cline{2-13}
      & \textbf{H0} & \textbf{H1} & \textbf{H0+H1}
      & \textbf{H0} & \textbf{H1} & \textbf{H0+H1}
      & \textbf{H0} & \textbf{H1} & \textbf{H0+H1}
      & \textbf{H0} & \textbf{H1} & \textbf{H0+H1} \\
    \hline\hline
    \textbf{Score} 
      & 0.16 & 2.89 & 1.37 
      & 2.76 & 0.41 & 0.72 
      & 0.58 & 1.00 & 0.21 
      & 0.10 & 12.0 & 12.0 \\
    \textbf{Resolution} 
      & 100 & 150 & 150 
      & 150 & 150 & 100 
      & 150 & 150 & 100 
      & (10,10) & (50,50) & (50,50) \\
    \textbf{Smoothing} 
      & 2 & 2 & 2 
      & 1 & 0 & 2 
      & 0 & 2 & 2 
      & 0 & 0 & 0 \\
    \textbf{Threshold} 
      & p10 & p10 & p10 
      & 0   & 0   & 0 
      & p10 & p10 & p10 
      & p10 & 0   & p10 \\
    \textbf{Normalization} 
      & $L^1$ & None & $L^1$ 
      & None  & $L^1$ & $L^1$ 
      & None  & $L^1$ & None 
      & $L^1$ & $L^1$ & $L^1$ \\
    \textbf{Landscapes} 
      &  &  &  
      & 5  & 1  & 1 
      &  &  &  
      &  &  &  \\
    \textbf{Bandwidth} 
      &  &  &  
      &  &  &  
      &  &  &  
      & 2   & 0.5 & 2 \\
    \textbf{Weights} 
      &  &  &  
      &  &  &  
      &  &  &  
      & const & const & const \\
    \hline\hline
    \end{tabular}%
    }
    \caption{Best hyperparameter grid settings (Table~\ref{grid}) for each topological representation and dimension.}
    \label{mejoresgrid}
\end{table}

 \begin{figure}[H]
         \centering
         \includegraphics[width=1\textwidth]{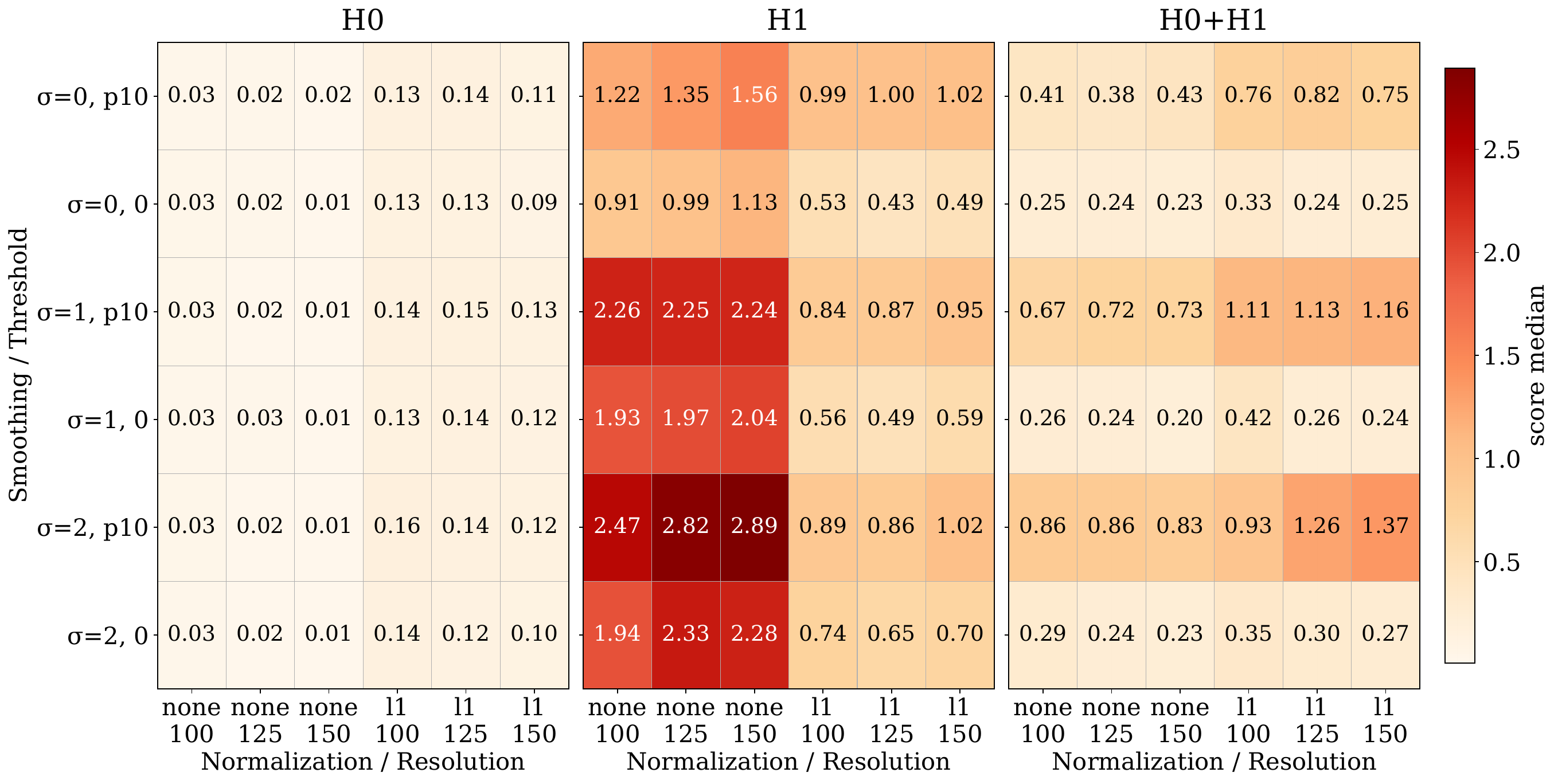}
             \caption{Parameter grid of median MWU scores for Betti curves across hyperparameter settings, stratified by homological dimension.}
         \label{fig:gridBC01}
     \end{figure}

 \begin{figure}[H]
         \centering
         \includegraphics[width=1\textwidth]{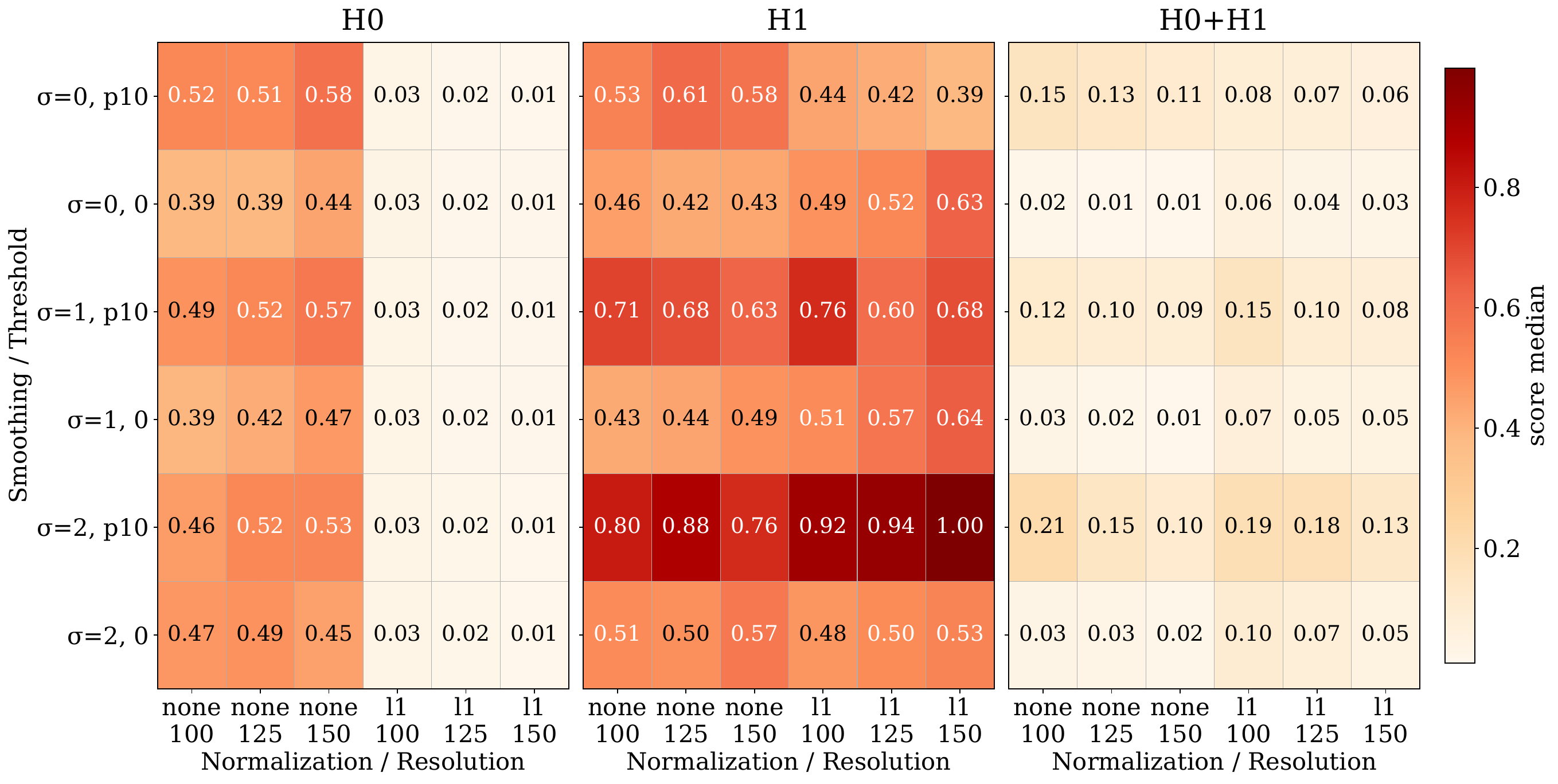}
             \caption{Parameter grid of median MWU scores for persistence silhouettes across hyperparameter settings, stratified by homological dimension.}
         \label{fig:gridBC02}
     \end{figure}

 \begin{figure}[H]
         \centering
         \includegraphics[width=1\textwidth]{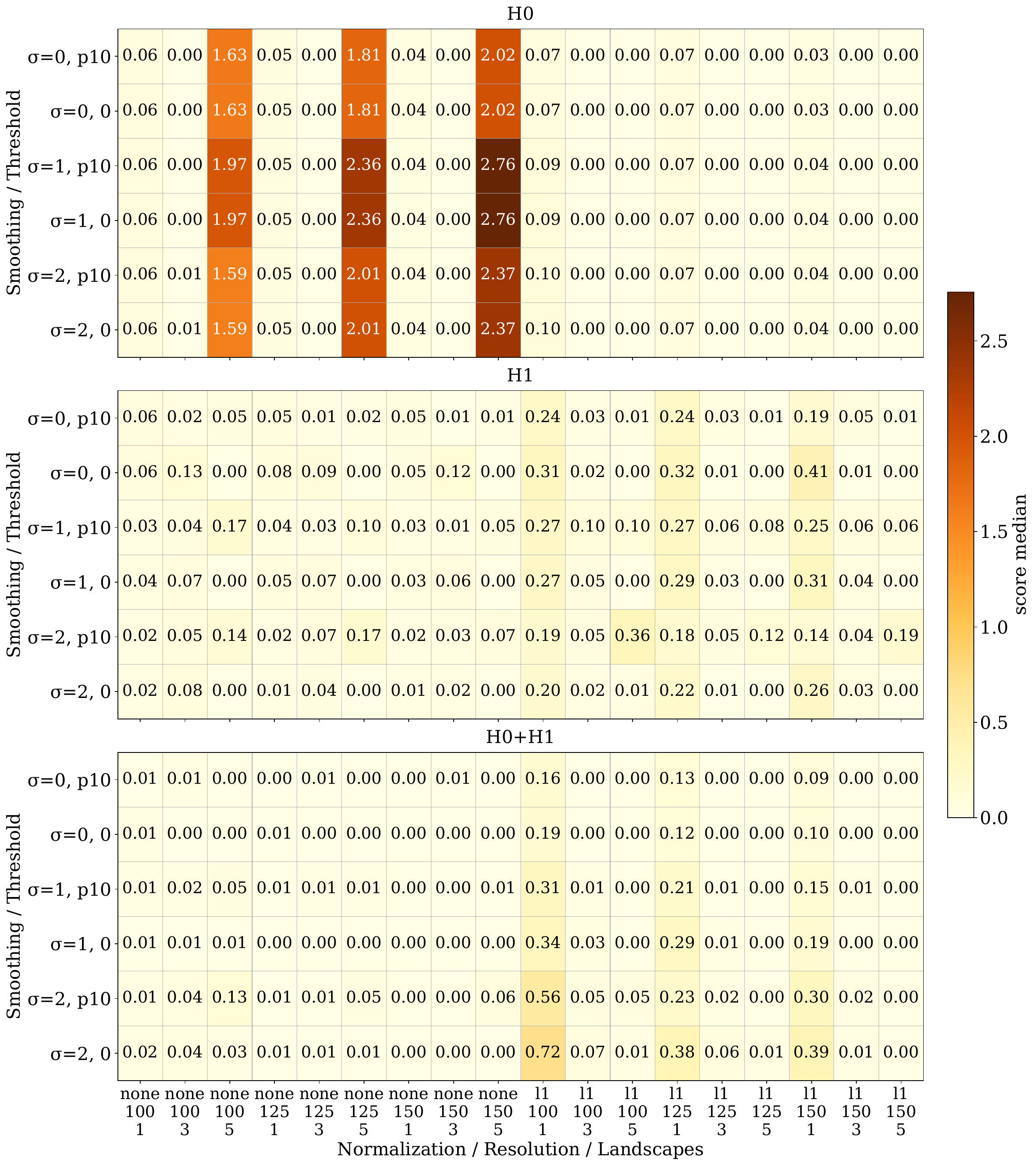}
             \caption{Parameter grid of median MWU scores for persistence landscapes across hyperparameter settings, stratified by homological dimension.}
         \label{fig:gridBC03}
     \end{figure}

\begin{figure}[H]
         \centering
         \includegraphics[width=1\textwidth]{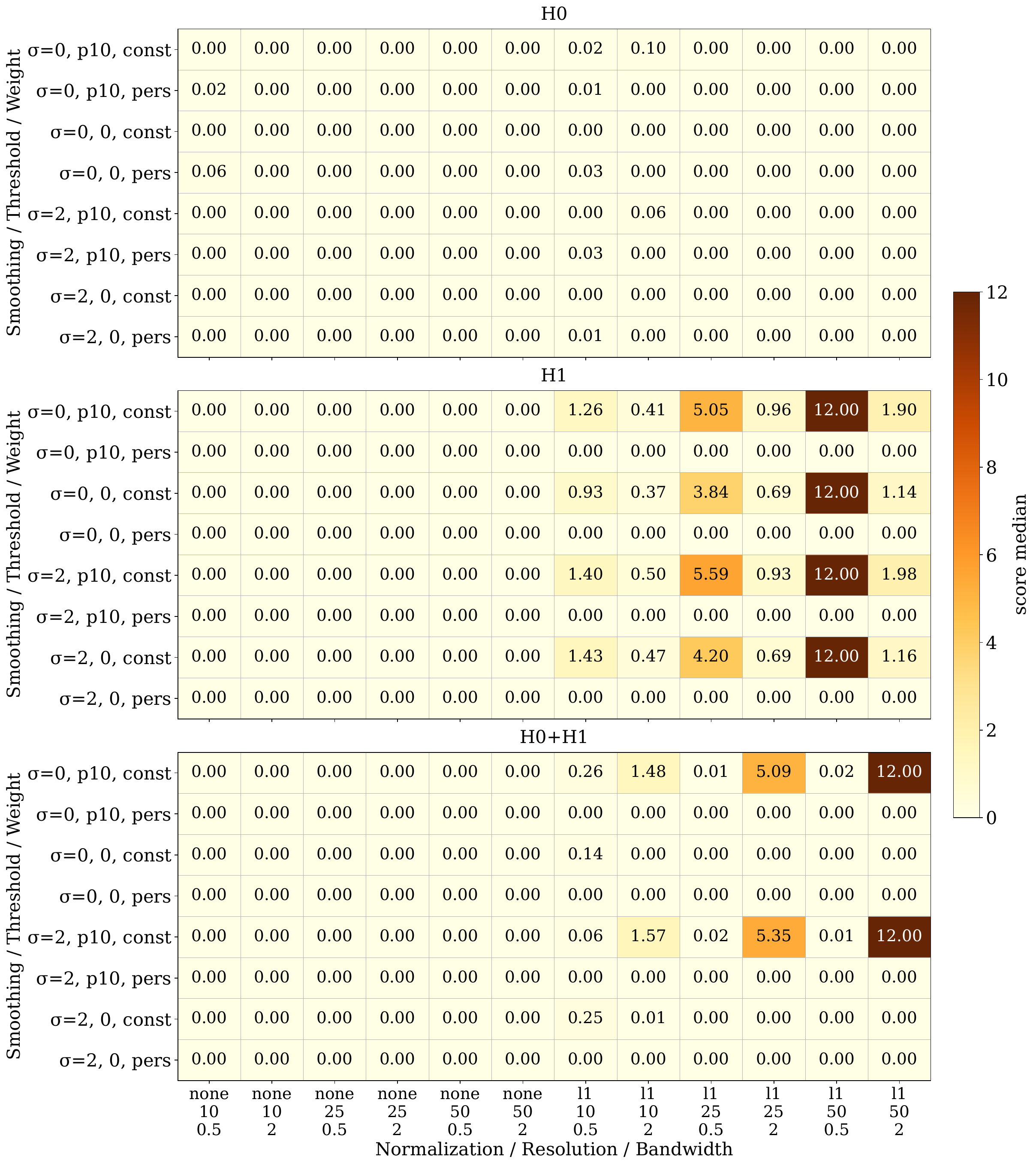}
             \caption{Parameter grid of median MWU scores for persistence images across hyperparameter settings, stratified by homological dimension.}
         \label{fig:gridBC04}
     \end{figure}

\begin{figure}[H]
         \centering
         \includegraphics[width=0.95\textwidth]{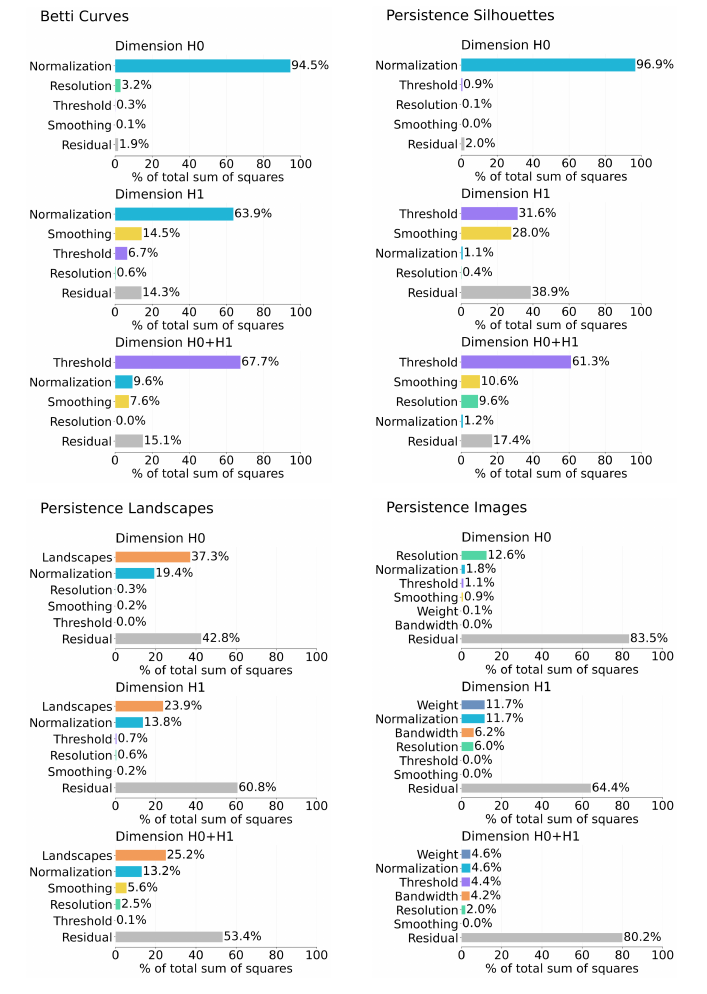}
             \caption{Decomposition of variance by TDA parameter across Betti curves, persistence silhouettes, persistence landscapes and persistence images, stratified by homological dimension.}
         \label{fig:anova_dims}
     \end{figure}

\subsection*{\(R_2\)}
\label{supp-r2}

\renewcommand{\arraystretch}{1.4}  
\begin{table}[H]
    \centering
    \small
    \renewcommand{\arraystretch}{1.3}
    \resizebox{\columnwidth}{!}{%
    \begin{tabular}{||l||c|c|c|c|c||}
    \hline\hline
        \textbf{Hyperparameter} 
        & \textbf{Random Forest} 
        & \textbf{SVM (linear)} 
        & \textbf{SVM (RBF)} 
        & \textbf{SVM (poly)} 
        & \textbf{SVM (sigmoid)} \\
    \hline\hline
        $n_{\text{estimators}}$ 
        & (50, 100, 200) 
        &  &  &  &  \\
        $\text{max\_depth}$ 
        & (None, 10) 
        &  &  &  &  \\
    \hline \hline
        $C$ 
        &  
        & (0.1, 10) 
        & (0.1, 10) 
        & (0.1, 10) 
        & (0.1, 10) \\
        $\gamma$ 
        &  
        &  
        & (0.01, 0.1) 
        & (0.01, 0.1) 
        & (0.01, 0.1) \\
        $\text{degree}$ 
        &  
        &  
        &  
        & (2, 3) 
        &  \\
        $\text{coef0}$ 
        &  
        &  
        &  
        & (0, 1) 
        & (0, 1) \\
    \hline\hline
    \end{tabular}%
    }
    \caption{Hyperparameter grid considered for the Random Forest and SVM classifiers. For each classifier, different combinations of parameters are evaluated in order to identify the settings that maximize the F2 score for classification between the NR and R groups.}
    \label{grid_clasificadores}

\end{table}

\begin{table}[H]
\centering
\footnotesize
\renewcommand{\arraystretch}{1.25}
\setlength{\tabcolsep}{4pt}

\begin{tabular}{|l|c|l|c|c|c|c|c|}
\hline
\textbf{Model} & \textbf{Oversamp.} & \textbf{Optimal Conf.} 
& \textbf{Dim} & \textbf{AUC} & \textbf{ACC} & \textbf{F2} & \textbf{CM} \\
\hline\hline

\multirow{6}{*}{\textbf{RF}} 

& \multirow{3}{*}{No}
& \begin{tabular}[c]{@{}l@{}}N\_est. = 50
\end{tabular}
& H0   & 0.66 & 0.79 & 0.20 & $\begin{pmatrix}72 & 4 \\ 16 & 4\end{pmatrix}$ \\ \cline{3-8}

&  
& \begin{tabular}[c]{@{}l@{}}N\_est. = 50
\end{tabular}
& H1   & 0.48 & 0.77 & 0.00 & $\begin{pmatrix}72 & 4 \\ 20 & 0\end{pmatrix}$ \\ \cline{3-8}

&  
& \begin{tabular}[c]{@{}l@{}}N\_est. = 50
\end{tabular}
& H0+H1 & 0.67 & 0.79 & 0.10 & $\begin{pmatrix}74 & 2 \\ 18 & 2\end{pmatrix}$ \\ \cline{2-8}

& \multirow{3}{*}{Yes}
& \begin{tabular}[c]{@{}l@{}}N\_est. = 200
\end{tabular}
& H0   & 0.70 & 0.79 & 0.35 & $\begin{pmatrix}69 & 7 \\ 13 & 7\end{pmatrix}$ \\ \cline{3-8}

&  
& \begin{tabular}[c]{@{}l@{}}N\_est. = 200
\end{tabular}
& H1   & 0.55 & 0.70 & 0.11 & $\begin{pmatrix}65 & 11 \\ 18 & 2\end{pmatrix}$ \\ \cline{3-8}

&  
& \begin{tabular}[c]{@{}l@{}}N\_est. = 100
\end{tabular}
& H0+H1 & 0.57 & 0.76 & 0.15 & $\begin{pmatrix}70 & 6 \\ 17 & 3\end{pmatrix}$ \\ 
\hline\hline

\multirow{6}{*}{\textbf{SVM}} 

& \multirow{3}{*}{No}
& \begin{tabular}[c]{@{}l@{}}Kernel = RBF\\ C = 0.1\\ $\gamma$ = 0.01\end{tabular}
& H0   & 0.67 & 0.54 & 0.54 & $\begin{pmatrix}38 & 38 \\ 6 & 14\end{pmatrix}$ \\ \cline{3-8}

&  
& \begin{tabular}[c]{@{}l@{}}Kernel = Sigmoid\\ C = 0.1\\ $\gamma$ = 0.01\\ coef0 = 1\end{tabular}
& H1   & 0.47 & 0.26 & 0.44 & $\begin{pmatrix}10 & 66 \\ 5 & 15\end{pmatrix}$ \\ \cline{3-8}

&  
& \begin{tabular}[c]{@{}l@{}}Kernel = RBF\\ C = 0.1\\ $\gamma$ = 0.01\end{tabular}
& H0+H1 & 0.61 & 0.35 & 0.48 & $\begin{pmatrix}19 & 57 \\ 5 & 15\end{pmatrix}$ \\ \cline{2-8}

& \multirow{3}{*}{Yes}
& \begin{tabular}[c]{@{}l@{}}Kernel = Linear\\ C = 0.1\end{tabular}
& H0   & 0.72 & 0.72 & 0.49 & $\begin{pmatrix}58 & 18 \\ 9 & 11\end{pmatrix}$ \\ \cline{3-8}

&  
& \begin{tabular}[c]{@{}l@{}}Kernel = Sigmoid\\ C = 10\\ $\gamma$ = 0.1\\ coef0 = 1\end{tabular}
& H1   & 0.51 & 0.46 & 0.41 & $\begin{pmatrix}33 & 43 \\ 9 & 11\end{pmatrix}$ \\ \cline{3-8}

&  
& \begin{tabular}[c]{@{}l@{}}Kernel = Sigmoid\\ C = 10\\ $\gamma$ = 0.1\\ coef0 = 1\end{tabular}
& H0+H1 & 0.47 & 0.45 & 0.47 & $\begin{pmatrix}30 & 46 \\ 7 & 13\end{pmatrix}$ \\ 
\hline

\end{tabular}

\caption{Classification results for Betti Curves generated using the parameters from Table~\ref{mejoresgrid}, the best Random Forest and SVM configurations, with and without oversampling.}
\label{tab:bc_rf_svm}
\end{table}

\begin{table}[H]
\centering
\footnotesize
\renewcommand{\arraystretch}{1}
\setlength{\tabcolsep}{4pt}

\begin{tabular}{|l|c|l|c|c|c|c|c|}
\hline
\textbf{Model} & \textbf{Oversamp.} & \textbf{Configuration} 
& \textbf{Dim} & \textbf{AUC} & \textbf{ACC} & \textbf{F2} & \textbf{CM} \\
\hline\hline

\multirow{6}{*}{\textbf{RF}} 

& \multirow{3}{*}{Without}
& \begin{tabular}[c]{@{}l@{}}N\_est. = 50
\end{tabular}
& H0   & 0.68 & 0.82 & 0.27 & $\begin{pmatrix}74 & 2 \\ 15 & 5\end{pmatrix}$ \\ \cline{3-8}

&  
& \begin{tabular}[c]{@{}l@{}}N\_est. = 50
\end{tabular}
& H1   & 0.50 & 0.78 & 0.00 & $\begin{pmatrix}75 & 1 \\ 20 & 0\end{pmatrix}$ \\ \cline{3-8}

&  
& \begin{tabular}[c]{@{}l@{}}N\_est. = 100
\end{tabular}
& H0+H1 & 0.66 & 0.79 & 0.11 & $\begin{pmatrix}74 & 2 \\ 18 & 2\end{pmatrix}$ \\ \cline{2-8}

& \multirow{3}{*}{With}
& \begin{tabular}[c]{@{}l@{}}N\_est. = 50
\end{tabular}
& H0   & 0.68 & 0.79 & 0.31 & $\begin{pmatrix}70 & 6 \\ 14 & 6\end{pmatrix}$ \\ \cline{3-8}

&  
& \begin{tabular}[c]{@{}l@{}}N\_est. = 100
\end{tabular}
& H1   & 0.51 & 0.71 & 0.21 & $\begin{pmatrix}64 & 12 \\ 16 & 4\end{pmatrix}$ \\ \cline{3-8}

&  
& \begin{tabular}[c]{@{}l@{}}N\_est. = 200
\end{tabular}
& H0+H1 & 0.66 & 0.78 & 0.21 & $\begin{pmatrix}71 & 5 \\ 16 & 4\end{pmatrix}$ \\ 
\hline\hline

\multirow{6}{*}{\textbf{SVM}} 

& \multirow{3}{*}{Without}
& \begin{tabular}[c]{@{}l@{}}Kernel = RBF\\ C = 0.1\\ $\gamma$ = 0.01
\end{tabular}
& H0   & 0.70 & 0.51 & 0.55 & $\begin{pmatrix}34 & 42 \\ 5 & 15\end{pmatrix}$ \\ \cline{3-8}

&  
& \begin{tabular}[c]{@{}l@{}}Kernel = Sigmoid\\ C = 10\\ $\gamma$ = 0.1
\\ coef0 = 0\end{tabular}
& H1   & 0.53 & 0.50 & 0.51 & $\begin{pmatrix}34 & 42 \\ 6 & 14\end{pmatrix}$ \\ \cline{3-8}

&  
& \begin{tabular}[c]{@{}l@{}}Kernel = Poly\\ C = 0.1\\ $\gamma$ = 0.01\\ degree = 3\\ coef0 = 1\end{tabular}
& H0+H1 & 0.64 & 0.76 & 0.48 & $\begin{pmatrix}63 & 13 \\ 10 & 10\end{pmatrix}$ \\ \cline{2-8}

& \multirow{3}{*}{With}
& \begin{tabular}[c]{@{}l@{}}Kernel = Poly\\ C = 10\\ $\gamma$ = 0.01\\ degree = 3
\end{tabular}
& H0   & 0.64 & 0.73 & 0.42 & $\begin{pmatrix}61 & 15 \\ 11 & 9\end{pmatrix}$ \\ \cline{3-8}

&  
& \begin{tabular}[c]{@{}l@{}}Kernel = Poly\\ C = 0.1\\ $\gamma$ = 0.01\\ degree = 2\\ coef0 = 1\end{tabular}
& H1   & 0.71 & 0.65& 0.49 & $\begin{pmatrix}50 & 26 \\ 8 & 12\end{pmatrix}$ \\ \cline{3-8}

&  
& \begin{tabular}[c]{@{}l@{}}Kernel = Sigmoid\\ C = 10\\ $\gamma$ = 0.01
\\ coef0 = 0\end{tabular}
& H0+H1 & 0.51 & 0.45 & 0.41 & $\begin{pmatrix}32 & 44 \\ 9 & 11\end{pmatrix}$ \\ 
\hline

\end{tabular}

\caption{Classification results for Persistence Silhouettes: best Random Forest and SVM configurations, with and without oversampling, for the dimensions $H_0$, $H_1$ and $H_0{+}H_1$.}
\label{tab:rf_svm_Silhouettes}
\end{table}

\begin{table}[H]
\centering
\footnotesize
\renewcommand{\arraystretch}{1}
\setlength{\tabcolsep}{4pt}

\begin{tabular}{|l|c|l|c|c|c|c|c|}
\hline
\textbf{Model} & \textbf{Oversamp.} & \textbf{Configuration} 
& \textbf{Dim} & \textbf{AUC} & \textbf{ACC} & \textbf{F2} & \textbf{CM} \\
\hline\hline

\multirow{6}{*}{\textbf{RF}} 

& \multirow{3}{*}{Without}
& \begin{tabular}[c]{@{}l@{}}N\_est. = 50
\end{tabular}
& H0   & 0.64 & 0.81 & 0.13 & $\begin{pmatrix}76 & 0 \\ 18 & 2\end{pmatrix}$ \\ \cline{3-8}

&  
& \begin{tabular}[c]{@{}l@{}}N\_est. = 50
\end{tabular}
& H1   & 0.54 & 0.77 & 0.00 & $\begin{pmatrix}74 & 2 \\ 20 & 0\end{pmatrix}$ \\ \cline{3-8}

&  
& \begin{tabular}[c]{@{}l@{}}N\_est. = 50
\end{tabular}
& H0+H1 & 0.61 & 0.79 & 0.00 & $\begin{pmatrix}76 & 0 \\ 20 & 0\end{pmatrix}$ \\ \cline{2-8}

& \multirow{3}{*}{With}
& \begin{tabular}[c]{@{}l@{}}N\_est. = 50
\end{tabular}
& H0   & 0.58 & 0.75 & 0.21 & $\begin{pmatrix}68 & 8 \\ 16 & 4\end{pmatrix}$ \\ \cline{3-8}

&  
& \begin{tabular}[c]{@{}l@{}}N\_est. = 100
\end{tabular}
& H1   & 0.56 & 0.78 & 0.13 & $\begin{pmatrix}73 & 3 \\ 18 & 2\end{pmatrix}$ \\ \cline{3-8}

&  
& \begin{tabular}[c]{@{}l@{}}N\_est. = 200
\end{tabular}
& H0+H1 & 0.56 & 0.77 & 0.17 & $\begin{pmatrix}71 & 5 \\ 17 & 3\end{pmatrix}$ \\ 
\hline\hline

\multirow{6}{*}{\textbf{SVM}} 

& \multirow{3}{*}{Without}
& \begin{tabular}[c]{@{}l@{}}Kernel = Sigmoid\\ C = 10\\ $\gamma$ = 0.01
\\ coef0 = 1\end{tabular}
& H0   & 0.59 & 0.57 & 0.53 & $\begin{pmatrix}41 & 35 \\ 6 & 14\end{pmatrix}$ \\ \cline{3-8}

&  
& \begin{tabular}[c]{@{}l@{}}Kernel = Sigmoid\\ C = 10\\ $\gamma$ = 0.01
\\ coef0 = 1\end{tabular}
& H1   & 0.59 & 0.55 & 0.47 & $\begin{pmatrix}41 & 35 \\ 8 & 12\end{pmatrix}$ \\ \cline{3-8}

&  
& \begin{tabular}[c]{@{}l@{}}Kernel = Poly\\ C = 0.1\\ $\gamma$ = 0.01\\ degree = 2\\ coef0 = 0\end{tabular}
& H0+H1 & 0.55 & 0.25 & 0.58 & $\begin{pmatrix}4 & 72 \\ 0 & 20\end{pmatrix}$ \\ \cline{2-8}

& \multirow{3}{*}{With}
& \begin{tabular}[c]{@{}l@{}}Kernel = Sigmoid\\ C = 0.1\\ $\gamma$ = 0.01
\\ coef0 = 0\end{tabular}
& H0   & 0.63 & 0.53 & 0.53 & $\begin{pmatrix}37 & 39 \\ 6 & 14\end{pmatrix}$ \\ \cline{3-8}

&  
& \begin{tabular}[c]{@{}l@{}}Kernel = Sigmoid\\ C = 10\\ $\gamma$ = 0.01
\\ coef0 = 1\end{tabular}
& H1   & 0.53 & 0.44 & 0.44 & $\begin{pmatrix}30 & 46 \\ 8 & 12\end{pmatrix}$ \\ \cline{3-8}

&  
& \begin{tabular}[c]{@{}l@{}}Kernel = Poly\\ C = 0.1\\ $\gamma$ = 0.01\\ degree = 2\\ coef0 = 0\end{tabular}
& H0+H1 & 0.59 & 0.30 & 0.51 & $\begin{pmatrix}12 & 64 \\ 3 & 17\end{pmatrix}$ \\ 
\hline

\end{tabular}

\caption{Classification results for Persistence Landscapes: best Random Forest and SVM configurations, with and without oversampling, for the dimensions $H_0$, $H_1$ and $H_0{+}H_1$.}
\label{tab:rf_svm_Landscapes}
\end{table}

\begin{table}[H]
\centering
\footnotesize
\renewcommand{\arraystretch}{1}
\setlength{\tabcolsep}{4pt}

\begin{tabular}{|l|c|l|c|c|c|c|c|}
\hline
\textbf{Model} & \textbf{Oversamp.} & \textbf{Configuration} 
& \textbf{Dim} & \textbf{AUC} & \textbf{ACC} & \textbf{F2} & \textbf{CM} \\
\hline\hline

\multirow{6}{*}{\textbf{RF}} 

& \multirow{3}{*}{Without}
& \begin{tabular}[c]{@{}l@{}}N\_est. = 50
\end{tabular}
& H0   & 0.42 & 071 & 0.21 & $\begin{pmatrix}64 & 12 \\ 16 & 4\end{pmatrix}$ \\ \cline{3-8}

&  
& \begin{tabular}[c]{@{}l@{}}N\_est. = 50
\end{tabular}
& H1   & 0.59 & 0.70 & 0.05 & $\begin{pmatrix}66 & 10 \\ 19 & 1\end{pmatrix}$ \\ \cline{3-8}

&  
& \begin{tabular}[c]{@{}l@{}}N\_est. = 50
\end{tabular}
& H0+H1 & 0.55 & 0.77 & 0.22 & $\begin{pmatrix}70 & 6 \\ 16 & 4\end{pmatrix}$ \\ \cline{2-8}

& \multirow{3}{*}{With}
& \begin{tabular}[c]{@{}l@{}}N\_est. = 200
\end{tabular}
& H0   & 0.40 & 0.65 & 0.24 & $\begin{pmatrix}57 & 19 \\ 15 & 5\end{pmatrix}$ \\ \cline{3-8}

&  
& \begin{tabular}[c]{@{}l@{}}N\_est. = 200
\end{tabular}
& H1   & 0.59 & 0.71 & 0.39 & $\begin{pmatrix}60 & 16 \\ 12 & 8\end{pmatrix}$ \\ \cline{3-8}

&  
& \begin{tabular}[c]{@{}l@{}}N\_est. = 200
\end{tabular}
& H0+H1 & 0.54 & 0.75 & 0.33 & $\begin{pmatrix}70 & 6 \\ 16 & 4\end{pmatrix}$ \\ 
\hline\hline

\multirow{6}{*}{\textbf{SVM}} 

& \multirow{3}{*}{Without}
& \begin{tabular}[c]{@{}l@{}}Kernel = RBF\\ C = 0.1\\ $\gamma$ = 0.01
\end{tabular}
& H0   & 0.54 & 0.31 & 0.41 & $\begin{pmatrix}17 & 59 \\ 7 & 13\end{pmatrix}$ \\ \cline{3-8}

&  
& \begin{tabular}[c]{@{}l@{}}Kernel = Sigmoid\\ C = 0.1\\ $\gamma$ = 0.01
\\ coef0 = 1\end{tabular}
& H1   & 0.61 & 0.32 & 0.46 & $\begin{pmatrix}15 & 61 \\ 5 & 15\end{pmatrix}$ \\ \cline{3-8}

&  
& \begin{tabular}[c]{@{}l@{}}Kernel = Sigmoid\\ C = 10\\ $\gamma$ = 0.1
\\ coef0 = 0\end{tabular}
& H0+H1 & 0.67 & 0.59 & 0.45 & $\begin{pmatrix}46 & 30 \\ 9 & 11\end{pmatrix}$ \\ \cline{2-8}

& \multirow{3}{*}{With}
& \begin{tabular}[c]{@{}l@{}}Kernel = Linear\\ C = 0.1\\ $\gamma$ = None
\end{tabular}
& H0   & 0.51 & 0.57 & 0.38 & $\begin{pmatrix}46 & 30 \\ 11 & 9\end{pmatrix}$ \\ \cline{3-8}

&  
& \begin{tabular}[c]{@{}l@{}}Kernel = Sigmoid\\ C = 10\\ $\gamma$ = 0.1
\\ coef0 = 1\end{tabular}
& H1   & 0.57 & 0.55 & 0.50 & $\begin{pmatrix}40 & 36 \\ 7 & 13\end{pmatrix}$ \\ \cline{3-8}

&  
& \begin{tabular}[c]{@{}l@{}}Kernel = Sigmoid\\ C = 10\\ $\gamma$ = 0.01
\\ coef0 = 0\end{tabular}
& H0+H1 & 0.54 & 0.58 & 0.45 & $\begin{pmatrix}45 & 31 \\ 9 & 11\end{pmatrix}$ \\ 
\hline

\end{tabular}

\caption{Classification results for Persistence Images: best Random Forest and SVM configurations, with and without oversampling, for the dimensions $H_0$, $H_1$ and $H_0{+}H_1$.}
\label{tab:rf_svm_Images}
\end{table}

\begin{figure}[H]
         \centering
         \includegraphics[width=0.95\textwidth]{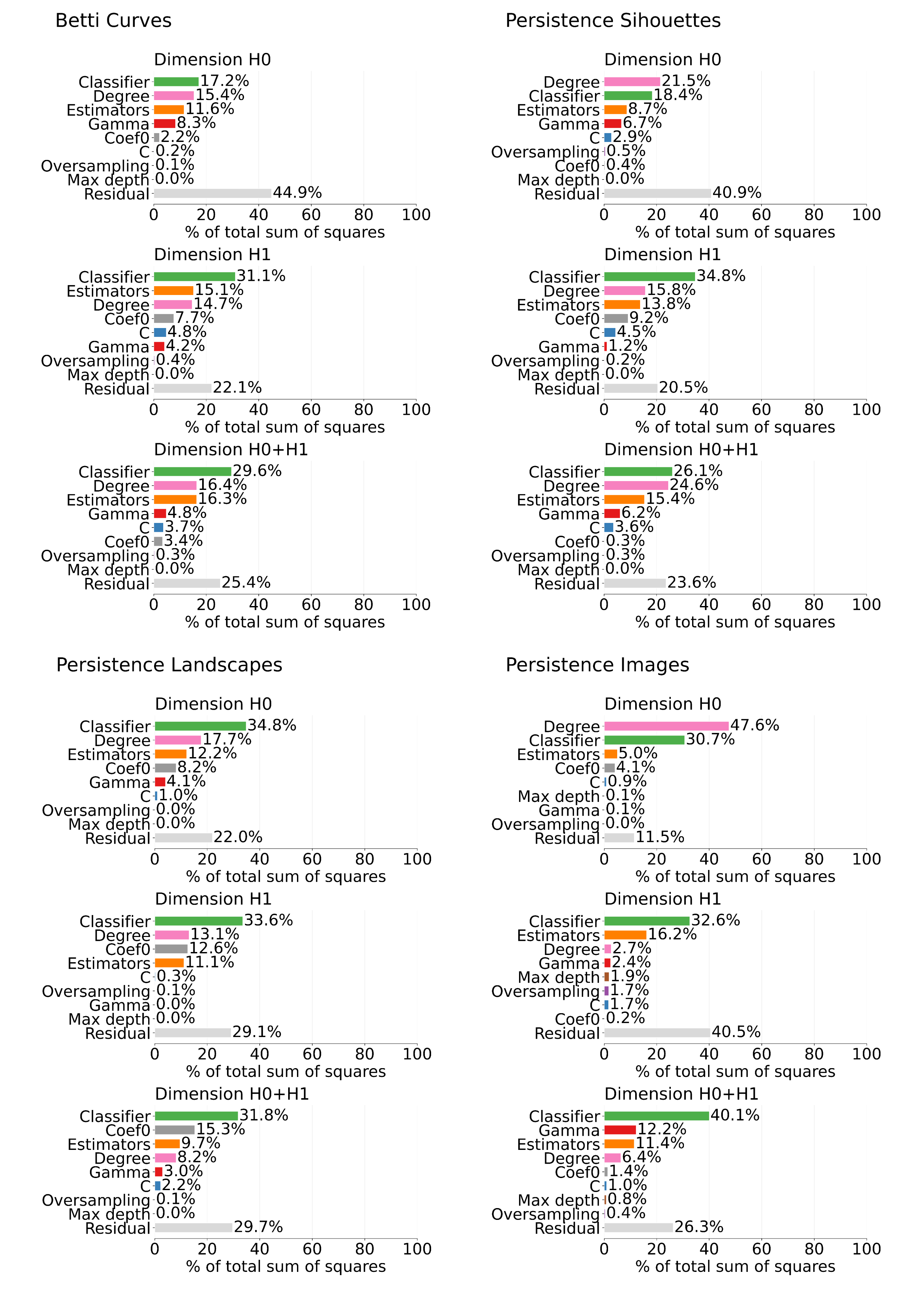}
             \caption{Decomposition of variance by TDA parameter across Betti Curves, Persistence Silhouettes, Persistence Landscapes and Persistence Images, stratified by homological dimension.}
         \label{fig:anova_dims}
     \end{figure}



\clearpage

  \bibliographystyle{elsarticle-harv} 
  \bibliography{cas-refs}



\end{document}